\documentclass[hidelinks,onefignum,onetabnum]{siamart250211}

\newcommand{\tp}{}

\usepackage{lipsum}
\usepackage{amsfonts}
\usepackage{graphicx}
\usepackage{epstopdf}
\usepackage{algorithmic}
\ifpdf
  \DeclareGraphicsExtensions{.eps,.pdf,.png,.jpg}
\else
  \DeclareGraphicsExtensions{.eps}
\fi

\newsiamremark{remark}{Remark}
\newsiamremark{hypothesis}{Hypothesis}
\crefname{hypothesis}{Hypothesis}{Hypotheses}
\newsiamthm{claim}{Claim}
\newsiamremark{fact}{Fact}
\crefname{fact}{Fact}{Facts}

\newsiamthm{assumption}{Assumption}
\crefname{assumption}{Assumption}{Assumptions}

\headers{Split-Merge for the Dominant Eigenvalue Problem}{X. Liu, M. Song, and Y. Xia}

\title{Split-Merge: A Difference-based Approach for Dominant Eigenvalue Problem\thanks{\funding{This research was supported by the National Key Research and Development Program of China (grant 2021YFA1003303), the National Natural Science Foundation of China (grant 12501421),  the Major Program of the National Natural Science Foundation of China (grants 72192830, 72192831), and the 111 Project (grant B16009).}}}

\author{Xiaozhi Liu\thanks{LMIB of the Ministry of Education, School of Mathematical Sciences, Beihang University, Beijing 100191, People's Republic of China (\email{xzliu@buaa.edu.cn},\ \email{yxia@buaa.edu.cn}).}
	\and Mengmeng Song\thanks{National Frontiers Science Center for Industrial Intelligence and Systems Optimization, Northeastern University, Shenyang 110819, People's Republic of China (\email{songmengmeng@mail.neu.edu.cn}).}
	\and Yong Xia \footnotemark[2] \thanks{Corresponding author.}}

\usepackage{amsopn}

\usepackage{color}

\usepackage{colortbl}
\definecolor{bblue}{rgb}{0.855,0.933,0.98}

\definecolor{rred}{HTML}{DC143C}

\usepackage{multirow}

\usepackage{makecell}
\usepackage{booktabs}
\usepackage{subcaption}
\usepackage{rotating}
\usepackage{siunitx}
\definecolor{smblue}{RGB}{221,238,255}
\definecolor{eigsgray}{RGB}{232,232,232}
\definecolor{dimrow}{RGB}{242,242,242}

\ifpdf
\hypersetup{
	pdftitle={Split-Merge: A Difference-based Approach for Dominant Eigenvalue Problem},
	pdfauthor={Xiaozhi Liu, Mengmeng Song, and Yong Xia}
}
\fi

\begin{document}
	
	\maketitle
	
	\begin{abstract}
		\tp
		The computation of the dominant eigenpair for symmetric positive semidefinite matrices is fundamental in numerical optimization. This work shifts the paradigm from the classical Rayleigh quotient to an unconstrained difference formulation, whose global optimum recovers the dominant eigenpair. Within this framework, we prove that gradient descent with a constant step-size $\alpha \in (0, 1)$ converges almost surely to the global optimum at a local linear rate. This analysis thereby reinterprets the classical power method as the conservative special case $\alpha=1/2$ and rigorously establishes its asymptotic sub-optimality. To advance this first-order scheme, we propose the Split-Merge algorithm based on the majorization-minimization principle. After splitting the matrix, we introduce auxiliary vectors to effectively merge the decomposition factors, resulting in a matrix-free and parameter-free iteration that captures tighter curvature information. We establish that Split-Merge converges almost surely to a global minimizer, and show that the iteration exhibits a spectral peeling mechanism that suppresses the targeted eigenspace, potentially surpassing the static linear rate of power iterations. Numerical evaluations across synthetic and real-world datasets confirm that our method has scalable efficiency, achieving speed-ups exceeding $10\times$ over the power method, with performance comparable to subspace iterations.
	\end{abstract}
	
	
	
	\section{Introduction}
	
	Computing the dominant eigenpair of a symmetric positive semidefinite (PSD) matrix is a central task in numerical optimization and linear algebra. It serves as a core subroutine in a wide range of scientific and engineering applications, such as principal component analysis (PCA) \cite{pca_2022}, spectral clustering \cite{spectral_clustering_2010}, PageRank \cite{pagerank_1999}, and low-rank matrix approximations \cite{low_rank_2013}.
	
	In PCA, the goal is to extract the dominant eigenpair of the sample covariance matrix $\boldsymbol{A} = \frac{1}{N} \sum_{i=1}^{N} \boldsymbol{d}_i \boldsymbol{d}_i^T \in \mathbb{R}^{n\times n}$.
	This well-studied task is classically cast as the maximization of the Rayleigh quotient \cite{rayleigh1896theory}:
	\begin{equation}
		\max_{\boldsymbol{x}\in \mathbb{R}^n} \frac{\boldsymbol{x}^T \boldsymbol{A} \boldsymbol{x}}{\boldsymbol{x}^T\boldsymbol{x}},\quad \text{s.t.} \quad \boldsymbol{x}\neq \boldsymbol{0}.
		\label{eq:rq_problem}
	\end{equation}
	The global maximum of \cref{eq:rq_problem} yields the dominant eigenvalue of the symmetric PSD matrix $\boldsymbol{A}$, with the optimum attained at the associated dominant eigenvector.

	This quotient formulation serves as the optimization foundation for most existing eigenvalue solvers, collectively referred to as quotient-based methods.
	
	\paragraph{Power Method and its Variations}
	
	The classical power method \cite{pm_1929} is widely adopted due to its matrix-free nature, yet its convergence rate degrades severely in the presence of a small spectral gap. To accelerate this process, several advanced variants have been proposed. Bai et al. \cite{ppm} introduced parameterized power methods (PPM) by applying gradient descent to the Rayleigh quotient, while Xu et al. \cite{power_m} proposed the power method with momentum (Power+M) inspired by the heavy ball method \cite{heavy_ball_polyak}. A persistent challenge across these accelerated schemes is their dependence on exact spectral priors. Specifically, the optimal step-size in PPM and the optimal momentum parameter in Power+M both require explicit knowledge of the sub-dominant eigenvalues. Although dynamic heuristic variants \cite{dynamic_power_m,DMPower} attempt to estimate these parameters on the fly, their rigorous convergence analyses still depend on unavailable spectral assumptions, limiting their theoretical completeness and practical utility.
	
	\paragraph{Advanced Subspace Methods}
	
	{
	\tp
	More sophisticated approaches leverage subspace projections to reduce computational complexity, including the Lanczos method~\cite{lanczos}, the locally optimal block preconditioned conjugate gradient (LOBPCG) method~\cite{LOBPCG}, and the Jacobi-Davidson (JD) method~\cite{jd_2000}.
	The Lanczos method constructs a Krylov subspace to iteratively project a symmetric matrix $\boldsymbol{A}$ onto a low-dimensional tridiagonal form, enabling efficient eigenvalue computation. However, a major algorithmic limitation of this process in finite-precision arithmetic is the rapid loss of orthogonality among the basis vectors. Maintaining orthogonality requires storing the entire Krylov basis for explicit reorthogonalization, which incurs substantial memory costs as the subspace grows and ultimately necessitates restarting strategies.
	Conversely, LOBPCG avoids global tridiagonalization and long-term memory dependencies by operating within a preconditioned conjugate gradient framework. It achieves local optimality by applying the Rayleigh-Ritz procedure over a low-dimensional trial subspace spanned by the current approximations, the preconditioned residuals, and the previous search directions. While highly robust and memory-efficient, each iteration of LOBPCG inherently necessitates the explicit orthonormalization of the local basis vectors and the solution of a small dense generalized eigenvalue subproblem. Consequently, its overall efficiency heavily relies on the availability of an effective preconditioner to keep the iteration count low.
	The JD method takes a different approach by combining the Rayleigh-Ritz procedure with flexible subspace expansion. It recasts the eigenvalue problem as a nonlinear system solved via approximate Newton iterations, but each step requires solving a linear system with a varying coefficient matrix, incurring high computational cost.
	}
	
	In contrast, this work explores the eigenvalue problem by minimizing the Auchmuty difference \cite{auchmuty1983duality}:
	\begin{equation}
		\min_{\boldsymbol{x} \in \mathbb{R}^n} \; f(\boldsymbol{x}) := \|\boldsymbol{x}\|^2 - (\boldsymbol{x}^T \boldsymbol{A} \boldsymbol{x})^{\frac{1}{2}}.
		\label{eq:dc_problem}
	\end{equation}
	Serving as a variational counterpart to the classical constrained quotient formulation \cref{eq:rq_problem} \cite{auchmuty1983duality,hiriart1985generalized}, this difference framework elegantly translates a fundamental numerical linear algebra task into an unconstrained optimization problem. Consequently, it has inspired a diverse class of difference-based solvers.
	
	\paragraph{Difference-based Methods}
	Early explorations by Mongeau and Torki \cite{mongeau2004computing} directly adapted classical steepest descent and Newton methods to this unconstrained formulation. To accelerate convergence without full Hessian computations, subsequent strategies integrated advanced step-size rules and quasi-Newton updates, yielding the Barzilai-Borwein-like algorithms of Gao et al. \cite{gao2015barzilai} and the modified L-BFGS scheme of Shi et al. \cite{shi2016limited}.
	However, these solvers predominantly optimize a smoothed variant of \cref{eq:dc_problem}, governed by the quartic objective $\|\boldsymbol{x}\|^4 - \boldsymbol{x}^T \boldsymbol{A} \boldsymbol{x}$. Despite their empirical utility, these approaches lack theoretical guarantees of achieving global optimality.
	
	{
	\tp
	In this paper, we investigate the intrinsic structure of the difference formulation in \cref{eq:dc_problem}. Our main contributions are summarized as follows:
	
	\begin{itemize}			
		\item \textit{First-Order Equivalence and Asymptotic Acceleration.} We show the equivalence between the classical power method and gradient descent applied to \cref{eq:dc_problem} with a constant step-size of $1/2$. Leveraging this first-order perspective, we prove that for any step-size $\alpha \in (0,1)$, gradient descent almost surely converges to a global minimizer with a local linear rate $\varrho(\alpha)$. This analysis rigorously quantifies the asymptotic sub-optimality of the classical power method: given leading eigenvalues $\lambda_1 > \lambda_2 > 0$, the rate $\varrho(\alpha)$ is unimodal on $(1/2, 1)$, attaining its unique minimum at $\alpha^* = 1/(2 - \lambda_2/\lambda_1)$.
		
		\item \textit{The Matrix-Free Split-Merge Algorithm.} We identify the power iteration as an instance of the majorization-minimization framework, bounded by an isotropic curvature surrogate and restricted to purely first-order geometry. To break this limitation, we implicitly split the PSD matrix to construct a tighter, curvature-aware local surrogate. By maximizing the single-step descent of this surrogate via a tractable relaxation, we extract an auxiliary vector that merges the factorization factors, yielding a simple, matrix-free, and parameter-free iteration.
		
		\item \textit{Global Optimality and Adaptive Spectral Peeling.} We rigorously establish convergence to the global minimizer for the Split-Merge algorithm. By characterizing the update as an adaptive spectral filter, we detail a spectral peeling mechanism that selectively suppresses localized eigenspaces. This targeted suppression overcomes the static linear rate of classical power iterations.
		
		\item \textit{Empirical Validation of Scalability and Efficiency.} Extensive numerical evaluations across synthetic and real-world datasets corroborate the theoretical advantages of the proposed framework. The Split-Merge algorithm achieves speed-ups exceeding $10\times$ over standard power iterations, delivering computational efficiency that rivals or surpasses sophisticated subspace methods, such as Lanczos and LOBPCG.
	\end{itemize}
	}

	The paper proceeds as follows. \Cref{sec:preliminary} establishes essential preliminary results regarding the difference formulation. \Cref{sec:pm_is_gd} recasts the classical power method within a first-order optimization framework, yielding a new convergence analysis and an asymptotic acceleration strategy. \Cref{sec:alg} develops the proposed Split-Merge algorithm within the majorization-minimization framework, and \cref{sec:convergence} rigorously derives its theoretical guarantees, including convergence to the global minimizer and the underlying acceleration mechanism. Numerical experiments demonstrating the efficiency and scalability of the approach are provided in \cref{sec:experiment}. Finally, \cref{sec:conclusion} concludes the paper and outlines avenues for future research.
	
	\noindent \textbf{Notation.} 
	Throughout, $\boldsymbol{A}$ denotes a matrix, $\boldsymbol{a}$ a vector, and $a$ a scalar, with $\boldsymbol{A}^T$ and $\boldsymbol{A}^{-1}$ representing the transpose and inverse of $\boldsymbol{A}$, respectively.  
	The relations $\boldsymbol{A} \succ 0$ and $\boldsymbol{A} \succeq 0$ denote positive definiteness and positive semidefiniteness, respectively.
	For a vector $\boldsymbol{a}$, $\operatorname{diag}(\boldsymbol{a})$ denotes the diagonal matrix whose diagonal elements are given by $\boldsymbol{a}$, $(\boldsymbol{a})_I$ represents the subvector indexed by the set $I$, and $(\boldsymbol{a})_i$ is its $i$-th element.
	Let $\|\cdot\|$ denote the Euclidean $\ell_2$-norm.
	The orthogonal projection of a point $\boldsymbol{x}_0$ onto a subspace $\mathcal{Q}$ is given by $\mathcal{P}_{\mathcal{Q}}(\boldsymbol{x}_0) = \mathop{\mathrm{argmin}}_{\boldsymbol{x} \in \mathcal{Q}} \|\boldsymbol{x} - \boldsymbol{x}_0\|$.
	
	\section{Preliminaries}\label{sec:preliminary}
	
	This paper addresses the computation of the dominant eigenpair of a symmetric PSD matrix $\boldsymbol{A}$ via the unconstrained optimization problem \cref{eq:dc_problem}. The PSD assumption entails no loss of generality. The extraction of extreme or interior eigenvalues from an arbitrary symmetric matrix reduces to finding the dominant eigenpair of a mapped PSD matrix $r(\boldsymbol{A})$. Such mappings include $r(\boldsymbol{A}) = \boldsymbol{A} + \eta \boldsymbol{I}$ for indefinite matrices, $r(\boldsymbol{A}) = \eta \boldsymbol{I} - \boldsymbol{A}$ to extract the smallest eigenvalue, {\tp and $r(\boldsymbol{A}) = \eta \boldsymbol{I} - (\boldsymbol{A} - \xi \boldsymbol{I})^2$ for targeting an interior eigenvalue close enough to $\xi$.}
	
	{
	\tp
	Throughout this paper, we consider a symmetric PSD matrix $\boldsymbol{A} \in \mathbb{R}^{n \times n}$ with rank $r$ satisfying $1 \le r \le n$. Its spectral decomposition is given by $\boldsymbol{A} = \boldsymbol{Q}\boldsymbol{\Lambda}\boldsymbol{Q}^T$, where $\boldsymbol{Q} = [\boldsymbol{q}_1, \dots, \boldsymbol{q}_n]$ is an orthogonal matrix and $\boldsymbol{\Lambda} = \operatorname{diag}(\lambda_1, \dots, \lambda_n)$. The spectrum is ordered such that the positive eigenvalues satisfy:
	\begin{equation*}
		\lambda_1 = \dots = \lambda_{m_1} > \lambda_{m_1+1} = \dots = \lambda_{m_2} > \lambda_{m_2+1} \ge \dots \ge \lambda_r > 0,
	\end{equation*}
	with the remaining $n-r$ eigenvalues being zero. Here, $m_1$ and $m_2 - m_1$ denote the respective algebraic multiplicities of the first two dominant distinct eigenvalues.
	}
	
	The theoretical foundation of the proposed framework relies on the analytic properties of the unconstrained objective $f$, adapted from Auchmuty \cite{auchmuty1983duality}.
	
	\begin{lemma} \label{lem:non_differentiable}
		The set of differentiable points of $f$ is given by
		\begin{equation}\label{eq:Theta}
			\Theta = \left\{ \boldsymbol{x} \in \mathbb{R}^n \mid \boldsymbol{A}\boldsymbol{x} \neq \boldsymbol{0} \right\}.
		\end{equation}
		At any $\boldsymbol{x} \in \Theta$, the gradient and Hessian of $f$ are respectively evaluated as:
		\begin{align}
			\nabla f(\boldsymbol{x}) &= 2 \boldsymbol{x} - \frac{\boldsymbol{A}\boldsymbol{x}}{(\boldsymbol{x}^T\boldsymbol{A}\boldsymbol{x})^{\frac{1}{2}}}, \label{eq:grad_f} \\
			\nabla^2 f(\boldsymbol{x}) &= 2 \boldsymbol{I} - \frac{\boldsymbol{A}}{(\boldsymbol{x}^T\boldsymbol{A}\boldsymbol{x})^{\frac{1}{2}}} + \frac{(\boldsymbol{A}\boldsymbol{x}) (\boldsymbol{A}\boldsymbol{x})^T}{(\boldsymbol{x}^T\boldsymbol{A}\boldsymbol{x})^{\frac{3}{2}}}. \label{eq:hesse_f}
		\end{align}
	\end{lemma}
	
	\begin{theorem}[{\tp\cite[Theorem 8.1]{auchmuty1983duality}}] 
	\label{thm:optimality}
	
		(i) Every stationary point $\boldsymbol{x}$ of $f$ is an eigenvector of $\boldsymbol{A}$. The corresponding eigenvalue is given by $\lambda(\boldsymbol{x}) = 2 (\boldsymbol{x}^T \boldsymbol{A} \boldsymbol{x})^{\frac{1}{2}}$, {\tp with the Euclidean norm satisfying $ \|\boldsymbol{x}\|=\sqrt{\lambda(\boldsymbol{x})}/2$ and the objective value achieving $f(\boldsymbol{x})=-{\lambda(\boldsymbol{x})}/{4}$}.
		
		(ii) The global minimizers of \cref{eq:dc_problem} are the dominant eigenvectors {\tp scaled to $\|\boldsymbol{x}\| = \sqrt{\lambda_1}/2$,} yielding the optimal value $-\lambda_1/4$.
		
		(iii) All second-order stationary points of \cref{eq:dc_problem} are global minimizers. Consequently, every local minimizer is also a global minimizer. Equivalently, all {\tp properly scaled} eigenvectors of $\boldsymbol{A}$ not associated with the dominant eigenvalue correspond to strict saddle points of $f$.
	\end{theorem}
	
	Consequently, the dominant eigenvector can be recovered by computing a global minimizer of $f$, whose landscape is devoid of non-strict saddle points. Furthermore, evaluating \cref{eq:hesse_f} reveals that the Hessian of $f$ exhibits a bounded positive curvature. 
	
	\begin{lemma} \label{lem:hessian_bound}
		For any $\boldsymbol{x} \in \Theta$, the Hessian of $f$ satisfies:
		\begin{equation}
			\nabla^2 f(\boldsymbol{x}) \preceq 2\boldsymbol{I}.
			\label{eq:hessian_pos_lip_2}
		\end{equation}
	\end{lemma}
	
	\begin{proof}
		By \cref{eq:hesse_f}, evaluating the quadratic form yields:
		$$\boldsymbol{w}^T \left( 2\boldsymbol{I} - \nabla^2 f(\boldsymbol{x}) \right) \boldsymbol{w} = \boldsymbol{w}^T \left(\frac{\boldsymbol{A}}{(\boldsymbol{x}^T\boldsymbol{A}\boldsymbol{x})^{\frac{1}{2}}} - \frac{(\boldsymbol{A}\boldsymbol{x}) (\boldsymbol{A}\boldsymbol{x})^T}{(\boldsymbol{x}^T\boldsymbol{A}\boldsymbol{x})^{\frac{3}{2}}}\right) \boldsymbol{w} \geq 0, \quad \forall \boldsymbol{w} \in \mathbb{R}^n, $$
		where the non-negativity is guaranteed by the Cauchy-Schwarz inequality, immediately yielding \cref{eq:hessian_pos_lip_2}.
	\end{proof}
	
	This curvature bound facilitates the integration of the power iteration into a first-order optimization framework.

	{
	\tp
	\section{The Power Method as a Gradient Descent Method}\label{sec:pm_is_gd}
	
	The power method \cite{pm_1929,pm_1998} stands as a cornerstone of numerical linear algebra, computing the dominant eigenvector of a matrix $\boldsymbol{A}$ via the iteration:
	\begin{equation}
		\boldsymbol{x}_{k+1} = \frac{\boldsymbol{A} \boldsymbol{x}_k}{\|\boldsymbol{A} \boldsymbol{x}_k\|}.
		\label{eq:pm_update}
	\end{equation}
	
	Applying a gradient descent step to the objective function $f$ in \cref{eq:dc_problem} with a constant step-size of $1/2$ yields:
	\begin{equation} 
		\label{eq:gd_update}
		\boldsymbol{x}_{k+1} = \boldsymbol{x}_{k} - \frac{1}{2} \nabla f(\boldsymbol{x}_k) = \boldsymbol{x}_{k} - \frac{1}{2} \left( 2\boldsymbol{x}_k - \frac{\boldsymbol{Ax}_k}{(\boldsymbol{x}_k^T \boldsymbol{A} \boldsymbol{x}_k)^{\frac{1}{2}}} \right) = \frac{\boldsymbol{Ax}_k}{2(\boldsymbol{x}_k^T \boldsymbol{A} \boldsymbol{x}_k)^{\frac{1}{2}}}.
	\end{equation}
	Up to positive scaling, \cref{eq:gd_update} perfectly recovers the power iteration \cref{eq:pm_update}.
	
	This equivalence facilitates a rigorous convergence analysis of the power method via first-order optimization theory. Specifically, the classical scheme is bottlenecked by a conservative step-size of $1/2$. In what follows, we establish almost sure convergence to a global minimizer and characterize the local linear rate for any gradient descent update $\boldsymbol{x}_{k+1} = \boldsymbol{x}_k - \alpha \nabla f(\boldsymbol{x}_k)$ with a step-size $\alpha \in (0, 1)$.
	
	\begin{theorem}[Asymptotic Stationarity and Subsequential Convergence] \label{thm:asymptotic_stationarity}
		Let $\{\boldsymbol{x}_k\}$ be the sequence generated by the gradient descent update applied to the objective $f$ in \cref{eq:dc_problem}, with a constant step-size $\alpha \in (0, 1)$. If the initial point $\boldsymbol{x}_0$ satisfies $f(\boldsymbol{x}_0) < 0$, then the entire trajectory $\{\boldsymbol{x}_k\}$ strictly remains within the differentiable domain $\Theta$ defined in \cref{eq:Theta}. Furthermore, the objective value monotonically decreases, satisfying the sufficient decrease condition:
		\begin{equation}
			f(\boldsymbol{x}_{k+1}) \leq f(\boldsymbol{x}_k) - \alpha \left(1 -\alpha \right) \|\nabla f(\boldsymbol{x}_k)\|^2, \quad \forall k \ge 0.
			\label{eq:sufficient_decrease}
		\end{equation}
		Consequently, the sequence $\{\boldsymbol{x}_k\}$ is bounded and its gradient vanishes asymptotically:
		\begin{equation}
			\lim_{k \to \infty} \|\nabla f(\boldsymbol{x}_k)\| = 0.
			\label{eq:grad_limit}
		\end{equation}
		Moreover, every accumulation point of $\{\boldsymbol{x}_k\}$ is a stationary point of $f$.
	\end{theorem}
	
	\begin{proof}
		We first establish that the entire trajectory $\{\boldsymbol{x}_k\}$ remains strictly in $\Theta$. Suppose otherwise, letting $\hat{\boldsymbol{x}} \notin \Theta$ denote the first exit point. Thus, for some $k_0 \ge 1$, $\hat{\boldsymbol{x}}$ lies on the segment connecting $\boldsymbol{x}_{k_0-1}$ and $\boldsymbol{x}_{k_0}$, while all strictly preceding points along the path lie in $\Theta$.
		
		For any $\boldsymbol{x}, \boldsymbol{y} \in \mathbb{R}^n$ with their connecting segment in $\Theta$, Taylor's theorem with the exact integral remainder gives:
		$$
		f(\boldsymbol{y}) = f(\boldsymbol{x}) + \nabla f(\boldsymbol{x})^T (\boldsymbol{y} - \boldsymbol{x}) + \int_0^1 (1 - t) (\boldsymbol{y} - \boldsymbol{x})^T \nabla^2 f(\boldsymbol{x} + t(\boldsymbol{y} - \boldsymbol{x})) (\boldsymbol{y} - \boldsymbol{x}) dt.
		$$
		Bounding the Hessian via \cref{eq:hessian_pos_lip_2} yields the standard quadratic upper bound:
		\begin{equation}
			f(\boldsymbol{y}) \leq f(\boldsymbol{x}) + \nabla f(\boldsymbol{x})^T (\boldsymbol{y} - \boldsymbol{x}) + \|\boldsymbol{y} - \boldsymbol{x}\|^2.
			\label{eq:upper_lip}
		\end{equation}
		Applying \cref{eq:upper_lip} to the gradient descent step $\boldsymbol{x}_{k+1} = \boldsymbol{x}_k - \alpha \nabla f(\boldsymbol{x}_k)$ for $0 \le k \le k_0-2$ yields the one-step descent:
		$$
		f(\boldsymbol{x}_{k+1}) \leq f(\boldsymbol{x}_k) - \alpha \left(1 -\alpha \right) \|\nabla f(\boldsymbol{x}_k)\|^2.
		$$
		Given $\alpha \in (0, 1)$, this ensures monotonic descent along the valid trajectory, enforcing
		\begin{equation}
			f(\boldsymbol{x}_{k_0-1}) \le f(\boldsymbol{x}_{k_0-2}) \le \dots \le f(\boldsymbol{x}_0) < 0.
			\label{eq:f_negative_chain}
		\end{equation}
		Furthermore, any point $\boldsymbol{x}$ strictly between $\boldsymbol{x}_{k_0-1}$ and $\hat{\boldsymbol{x}}$ constitutes a partial step $\boldsymbol{x} = \boldsymbol{x}_{k_0-1} - t \nabla f(\boldsymbol{x}_{k_0-1})$ for some $t \in (0, \alpha)$. As this open segment lies in $\Theta$, invoking \cref{eq:upper_lip} again gives $f(\boldsymbol{x}) \le f(\boldsymbol{x}_{k_0-1})$. 
		
		By the continuity of $f$ and the descent in \cref{eq:f_negative_chain}, the limit as $\boldsymbol{x}$ tends to $\hat{\boldsymbol{x}}$ yields:
		\begin{equation}
			f(\hat{\boldsymbol{x}}) \le f(\boldsymbol{x}_{k_0-1}) < 0.
			\label{eq:f_hat_negative}
		\end{equation}
		However, $\hat{\boldsymbol{x}} \notin \Theta$ implies $\boldsymbol{A}\hat{\boldsymbol{x}} = \boldsymbol{0}$, trivially yielding $f(\hat{\boldsymbol{x}}) \ge 0$. This contradicts \cref{eq:f_hat_negative}, ensuring the entire trajectory remains in $\Theta$.
		Consequently, the bound \cref{eq:upper_lip} holds globally along the trajectory, yielding \cref{eq:sufficient_decrease} for all $k \ge 0$.
		
		Summing \cref{eq:sufficient_decrease} over $k = 0, \dots, N$ yields
		$$
		\alpha(1 - \alpha) \sum_{k=0}^N \|\nabla f(\boldsymbol{x}_k)\|^2 \leq f(\boldsymbol{x}_0) - f(\boldsymbol{x}_{N+1}).
		$$
		Because $f$ is coercive, it admits a finite global lower bound $f^* = \inf_{\boldsymbol{x} \in \mathbb{R}^n} f(\boldsymbol{x}) > -\infty$. Taking the limit as $N \to \infty$ for $\alpha \in (0, 1)$ establishes:
		$$
		\sum_{k=0}^\infty \|\nabla f(\boldsymbol{x}_k)\|^2 \leq \frac{f(\boldsymbol{x}_0) - f^*}{\alpha(1 - \alpha)} < \infty.
		$$
		The convergence of this series trivially necessitates \cref{eq:grad_limit}.
		
		Furthermore, the monotonic descent of $\{f(\boldsymbol{x}_k)\}$ confines the trajectory $\{\boldsymbol{x}_k\}$ to the sublevel set $\mathcal{L} := \{\boldsymbol{x} \in \mathbb{R}^n \mid f(\boldsymbol{x}) \leq f(\boldsymbol{x}_0)\}$. The coercivity of $f$ ensures $\mathcal{L}$ is bounded, thereby establishing the boundedness of $\{\boldsymbol{x}_k\}$ and guaranteeing the existence of at least one accumulation point.
		
		Let $\boldsymbol{x}^*$ be an arbitrary accumulation point, with a subsequence $\{\boldsymbol{x}_{k_j}\}$ converging to $\boldsymbol{x}^*$. By the continuous differentiability of $f$ in $\Theta$, taking the limit along the subsequence yields
		$$
		\nabla f(\boldsymbol{x}^*) = \lim_{j \to \infty} \nabla f(\boldsymbol{x}_{k_j}) = \boldsymbol{0}.
		$$
		Thus, every accumulation point of $\{\boldsymbol{x}_k\}$ is a stationary point of $f$.
	\end{proof}
	
	Having established the asymptotic vanishing of the gradient in \cref{thm:asymptotic_stationarity}, the following theorem elevates this subsequential behavior to sequence convergence.
	
	\begin{theorem}[Global Convergence to a Stationary Point] \label{thm:stationary_convergence}
		For any initial point satisfying $f(\boldsymbol{x}_0) < 0$, the sequence $\{\boldsymbol{x}_k\}$ generated by the gradient descent update for \cref{eq:dc_problem} with a constant step-size $\alpha \in (0, 1)$ converges to a nonzero stationary point of $f$.
	\end{theorem}
	
	\begin{proof}
		By \cref{thm:asymptotic_stationarity}, the sequence $\{f(\boldsymbol{x}_k)\}$ is monotonically decreasing and bounded below, establishing its convergence. Consequently, all accumulation points of $\{\boldsymbol{x}_k\}$ share the same negative objective value. Invoking \cref{thm:asymptotic_stationarity} alongside \cref{thm:optimality} (i), these accumulation points must be eigenvectors of $\boldsymbol{A}$ associated with a common positive eigenvalue, denoted as $\lambda^* > 0$. Furthermore, every accumulation point possesses the uniform norm $\sqrt{\lambda^*}/2$, yielding
		\begin{equation}
			\lim_{k\to\infty} (\boldsymbol{x}_k^T \boldsymbol{A} \boldsymbol{x}_k)^{\frac{1}{2}} = {\lambda^*}/{2}.
			\label{eq:limit_quadratic}
		\end{equation}
		
		Let $I := \{i \mid \lambda_i = \lambda^*\}$. Since all accumulation points lie in the eigenspace associated with $\lambda^*$, their projections onto the orthogonal complement vanish asymptotically:
		\begin{equation}
			\lim_{k\to\infty} \boldsymbol{q}_i^T \boldsymbol{x}_k = 0, \quad \forall i \notin I.
			\label{eq:orthogonal_vanish}
		\end{equation}
		Applying the gradient descent update $\boldsymbol{x}_{k+1} = \boldsymbol{x}_k - \alpha \nabla f(\boldsymbol{x}_k)$ with the explicit gradient from \cref{eq:grad_f} yields the projection iteration on $I$:
		\begin{equation}
			(\boldsymbol{Q}^T \boldsymbol{x}_{k+1})_I = \nu_k (\boldsymbol{Q}^T \boldsymbol{x}_k)_I, \quad \text{where} \quad \nu_k = 1 - \alpha \left( 2 - {\lambda^*}/{(\boldsymbol{x}_k^T \boldsymbol{A} \boldsymbol{x}_k)^{\frac{1}{2}}} \right).
			\label{eq:collinear_update}
		\end{equation}
		This iteration dictates that $(\boldsymbol{Q}^T \boldsymbol{x}_k)_I$ remains collinear with the initial projection $(\boldsymbol{Q}^T \boldsymbol{x}_0)_I$ for all $k \ge 1$. Given that all accumulation points of $\{\boldsymbol{x}_k\}$ possess the identical norm $\sqrt{\lambda^*}/2$, the sequence of projections $\{(\boldsymbol{Q}^T \boldsymbol{x}_k)_I\}$ admits at most a single accumulation vector and its negation.
		
		By \cref{eq:limit_quadratic}, the scaling factor satisfies $\lim_{k\to\infty} \nu_k = 1$. Consequently, $\nu_k$ remains positive for all sufficiently large $k$, rigorously precluding any sign oscillation. Coupled with the orthogonal vanishing in \cref{eq:orthogonal_vanish}, this guarantees that the entire trajectory $\{\boldsymbol{x}_k\}$ converges to a single nonzero stationary point.
	\end{proof}
	
	While \cref{thm:stationary_convergence} guarantees convergence to a nonzero stationary point, securing global optimality necessitates excluding suboptimal eigenvectors. By \cref{thm:optimality} (iii), every such suboptimal critical point constitutes a strict saddle of $f$. Exploiting this instability, we now prove that the initializations drawn to any non-global limit form a set of Lebesgue measure zero, ultimately sealing the almost sure convergence to a global minimizer.
	
	\begin{theorem}[Almost Sure Convergence to a Global Minimizer] \label{thm:global_minimizer}
		Let the initial point $\boldsymbol{x}_0$ be drawn from the set $\{\boldsymbol{x} \in \mathbb{R}^n \mid f(\boldsymbol{x}) < 0\}$ according to an absolutely continuous probability measure. Then, the entire sequence $\{\boldsymbol{x}_k\}$ generated by the gradient descent update for \cref{eq:dc_problem} with a constant step-size $\alpha \in (0, 1)$ converges to a global minimizer of $f$ with probability one.
	\end{theorem}
	
	\begin{proof}
		To establish almost sure convergence, it suffices to prove that the set of initializations yielding convergence to any strict saddle has Lebesgue measure zero. By \cref{thm:stationary_convergence}, assume that $\{\boldsymbol{x}_k\}$ converges to an eigenvector $\boldsymbol{x}^*$ associated with an eigenvalue $\lambda^* < \lambda_1$.
		
		Analogous to \cref{eq:collinear_update}, we obtain
		\begin{equation}
			\boldsymbol{q}_1^T \boldsymbol{x}_{k+1} = \nu_k \boldsymbol{q}_1^T \boldsymbol{x}_k, \quad \text{where} \quad \nu_k = 1 - \alpha \left( 2 - {\lambda_1}/{(\boldsymbol{x}_k^T \boldsymbol{A} \boldsymbol{x}_k)^{\frac{1}{2}}} \right).
			\label{eq:q1_recurrence}
		\end{equation}
		By \cref{eq:limit_quadratic}, the asymptotic scaling factor satisfies
		$$
		\lim_{k\to\infty} \nu_k = 1 + 2\alpha \left({\lambda_1}/{\lambda^*} - 1 \right) > 1.
		$$
		Fix $\epsilon \in (0, 2\alpha(\lambda_1/\lambda^* - 1))$. There exists $K_1 \ge 0$ such that $\nu_k > 1 + \epsilon$ for all $k \ge K_1$.
		
		If $\boldsymbol{q}_1^T \boldsymbol{x}_k \neq 0$ for all $k \ge 0$, the projection in \cref{eq:q1_recurrence} expands exponentially:
		$$
		|\boldsymbol{q}_1^T \boldsymbol{x}_k| > (1 + \epsilon)^{k - K_1} |\boldsymbol{q}_1^T \boldsymbol{x}_{K_1}|, \quad \forall k > K_1.
		$$
		This forces $|\boldsymbol{q}_1^T \boldsymbol{x}_k|$ to diverge to infinity, strictly contradicting the asymptotic behavior $\lim_{k\to\infty} \boldsymbol{q}_1^T \boldsymbol{x}_k = 0$ derived analogously to \cref{eq:orthogonal_vanish}.
		
		Define the base set $Z_0 = \{\boldsymbol{x}_0 \in \mathbb{R}^n \mid \boldsymbol{q}_1^T \boldsymbol{x}_0 = 0\}$, and for $K \ge 1$:
		$$
		Z_K = \{\boldsymbol{x}_0 \in \mathbb{R}^n \mid \boldsymbol{q}_1^T \boldsymbol{x}_K = 0 \text{ and } \boldsymbol{q}_1^T \boldsymbol{x}_k \neq 0, \, \forall k < K\}.
		$$
		The total set of initializations converging to any suboptimal stationary point is contained in $Z = \bigcup_{K=0}^\infty Z_K$.
		
		The base set $Z_0$ constitutes a hyperplane, trivially possessing Lebesgue measure zero. For $K \ge 1$, taking $\boldsymbol{x}_0 \in Z_K$ implies $\boldsymbol{q}_1^T \boldsymbol{x}_{K-1} \neq 0$ and $\boldsymbol{q}_1^T \boldsymbol{x}_K = 0$. Invoking \cref{eq:q1_recurrence} immediately yields
		\begin{equation}
			1 - 2\alpha + {\alpha \lambda_1}/{(\boldsymbol{x}_{K-1}^T \boldsymbol{A} \boldsymbol{x}_{K-1})^{\frac{1}{2}}} = 0.
			\label{eq:nu_zero_condition}
		\end{equation}
		The denominator is positive since the entire trajectory remains in $\Theta$ by \cref{thm:asymptotic_stationarity}. If $\alpha \in (0, 1/2]$, \cref{eq:nu_zero_condition} admits no real solution, giving $Z_K = \emptyset$.
		Otherwise, recursively substituting the algebraic gradient descent update $\boldsymbol{x}_k = \boldsymbol{x}_{k-1} - \alpha \nabla f(\boldsymbol{x}_{k-1})$ expresses \cref{eq:nu_zero_condition} as a real algebraic equation in terms of $\boldsymbol{x}_0$. Evaluating this relation at the global minimizer $\sqrt{\lambda_1}\boldsymbol{q}_1/2$ yields $1 - 2\alpha + 2\alpha = 1 \neq 0$, certifying its non-triviality. Since the solution set of any non-trivial real algebraic equation has Lebesgue measure zero, each $Z_K$ possesses measure zero. The countable union $Z$ thus has measure zero, almost surely precluding convergence to any such suboptimal limit.
	\end{proof}
	
	\begin{remark}[Randomized Initialization] \label{rmk:random_init}
		A valid initialization for \cref{thm:global_minimizer} can be constructed via uniform sampling. Suppose $\bar{\boldsymbol{p}} \in \mathbb{R}^n$ is drawn uniformly from the unit sphere, ensuring $\boldsymbol{A}\bar{\boldsymbol{p}} \neq \mathbf{0}$ almost surely. Choosing a scalar $r$ uniformly from the interval $(0, (\bar{\boldsymbol{p}}^T \boldsymbol{A} \bar{\boldsymbol{p}})^{\frac{1}{2}})$ and setting $\boldsymbol{x}_0 = r\bar{\boldsymbol{p}}$ ensures that $\boldsymbol{x}_0$ admits an absolutely continuous probability measure while rigorously satisfying $f(\boldsymbol{x}_0) < 0$.
	\end{remark}
	
	\begin{remark}[Connection to and Distinction from Classical Saddle-Escaping Theory] \label{rmk:saddle_escape}
		While sharing the spirit of generic saddle-escaping theories \cite{lee2019first, schaeffer2020extending}, a direct application of these results is incompatible with our specific problem structure. Such theories necessitate smoothness over convex domains and abstract spectral non-degeneracy to support large step-sizes. These prerequisites directly conflict with our non-convex domain $\Theta$ in \cref{eq:Theta}. By instead exploiting the exact algebraic recurrence of the eigenspace projection, our self-contained proof bypasses these structural limitations, yielding a structure-exploiting guarantee that establishes almost sure convergence to a global minimizer across the entire interval $\alpha \in (0, 1)$.
	\end{remark}
	
	Having established convergence to the global minimizer, the following theorem characterizes the local linear rate.
	
	\begin{theorem}[Local Linear Convergence] \label{thm:local_convergence}
		Assume $\lambda_1 > \lambda_2$. For any step-size $\alpha \in (0, 1)$, the gradient descent sequence $\{\boldsymbol{x}_k\}$ exhibits local linear convergence to a global minimizer $\boldsymbol{x}^*$, satisfying
		$$ \limsup_{k \to \infty} \frac{\|\boldsymbol{x}_{k+1} - \boldsymbol{x}^*\|}{\|\boldsymbol{x}_k - \boldsymbol{x}^*\|} = \varrho(\alpha). $$
		The asymptotic rate $\varrho(\alpha)$ is given by
		\begin{equation}
			\varrho(\alpha) = \max \left\{ |1 - 2\alpha|, \left| 1 - 2\alpha \left( 1 - {\lambda_2}/{\lambda_1} \right) \right| \right\} < 1.
			\label{eq:rate_alpha}
		\end{equation}
	\end{theorem}
	
	\begin{proof}
		The asymptotic rate is governed by the spectral radius $\varrho(\alpha)$ of the Jacobian operator $\boldsymbol{T}_\alpha = \boldsymbol{I} - \alpha \nabla^2 f(\boldsymbol{x}^*)$. Evaluating the exact Hessian in \cref{eq:hesse_f} at the canonical minimizer $\boldsymbol{x}^* = \sqrt{\lambda_1}\boldsymbol{q}_1/2$ yields
		\begin{equation*}
			\nabla^2 f(\boldsymbol{x}^*) = 2 \boldsymbol{I} - \frac{2}{\lambda_1}\boldsymbol{A} + 2\boldsymbol{q}_1\boldsymbol{q}_1^T.
		\end{equation*}
		Given the spectral decomposition $\boldsymbol{A} = \sum_{i=1}^n \lambda_i \boldsymbol{q}_i \boldsymbol{q}_i^T$, the eigensystem of $\nabla^2 f(\boldsymbol{x}^*)$ aligns with $\boldsymbol{Q}$. We extract the exact Hessian eigenvalues $\xi_i$:
		$$ \xi_1 = 2, \quad \text{and} \quad \xi_i = 2 \left( 1 - {\lambda_i}/{\lambda_1} \right), \quad \forall i \ge 2. $$
		For any $\alpha \in (0, 1)$, the eigenvalues $1 - \alpha \xi_i$ of $\boldsymbol{T}_\alpha$ reside in $(-1, 1)$, ensuring $\varrho(\alpha)=\max _{1 \leq i \leq n}\left|1-\alpha \xi_i\right| < 1$. Since $\xi_i \in [\xi_2, \xi_1]$, the spectral radius depends solely on these extremal bounds, establishing \cref{eq:rate_alpha}.
	\end{proof}
	
	\begin{corollary}[The Classical Step-Size $\alpha=1/2$ is Asymptotically Suboptimal] \label{cor:suboptimal}
		For $\lambda_1 > \lambda_2 > 0$, the local linear rate $\varrho(\alpha)$ in \cref{eq:rate_alpha} is unimodal on the interval $(1/2, 1)$. While setting $\alpha = 1/2$ exactly recovers the standard power method rate $\varrho(1/2) = \lambda_2/\lambda_1$, the asymptotically optimal step-size is given by $\alpha^* = 1/(2 - \lambda_2/\lambda_1) \in (1/2, 1)$, achieving the superior rate:
		\begin{equation}
			\varrho(\alpha^*) = \frac{\lambda_2/\lambda_1}{2 - \lambda_2/\lambda_1} < \frac{\lambda_2}{\lambda_1}.
			\label{eq:rate_optimal}
		\end{equation}
	\end{corollary}
	
	\begin{proof}
		Evaluating \cref{eq:rate_alpha} at $\alpha = 1/2$ yields $\varrho(1/2) = \lambda_2/\lambda_1$. On the interval $(1/2, 1)$, the rate $\varrho(\alpha)$ is unimodal, governed by the intersection of a strictly decreasing and a strictly increasing linear bound. Equating these bounds $2\alpha - 1 = 1 - 2\alpha(1 - \lambda_2/\lambda_1)$ uniquely minimizes the rate, isolating the optimal step-size $\alpha^* = 1 / (2 - \lambda_2/\lambda_1)$.
		Since $\lambda_1 > \lambda_2 > 0$, it follows that $\alpha^* \in (1/2, 1)$. Substituting $\alpha^*$ back into \cref{eq:rate_alpha} establishes the superior asymptotic bound, yielding~\cref{eq:rate_optimal}.
	\end{proof}
	
	While the asymptotically optimal step-size $\alpha^*$ remains unavailable without prior spectral knowledge, our primary contribution is to explicitly demonstrate the suboptimality of the power method within the first-order framework. The enhanced descent trajectory governed by the optimum $\alpha^*$ is corroborated in \cref{fig:fvalueitergd}.
	
	\begin{figure}[t]
		\centering
		\includegraphics[width=0.7\linewidth]{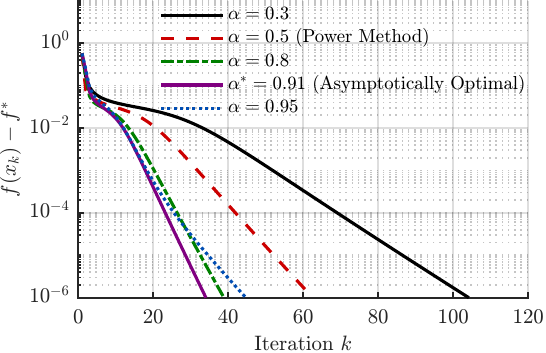}
		\caption{\tp \textbf{Asymptotic acceleration of the gradient descent update.} The performance is evaluated on a synthetic matrix ($n=1024$) with leading eigenvalues $\lambda_1 = 1$ and $\lambda_2 = 0.9$. The vertical axis reports the optimality gap $f(\boldsymbol{x}_k) - f^*$, where the optimal value is $f^* = -{\lambda_1}/{4}$.}
		\label{fig:fvalueitergd}
	\end{figure}
	}
	
	\section{The Split-Merge Algorithm} \label{sec:alg}
		
	{
	\tp
	While the gradient descent equivalence in \cref{sec:pm_is_gd} guarantees convergence to the global minimizer and exposes the conservativeness of the power method step-size, the iterative update remains confined to purely first-order geometry. We subsequently recast the algorithmic formulation within the majorization-minimization framework \cite{mm_2016}.
	
	Observe that the objective \cref{eq:dc_problem} admits a natural difference of convex (DC) decomposition: $f(\boldsymbol{x}) = g(\boldsymbol{x}) - h(\boldsymbol{x})$, where $g(\boldsymbol{x}) = \|\boldsymbol{x}\|^2$ and $h(\boldsymbol{x}) = (\boldsymbol{x}^T\boldsymbol{A}\boldsymbol{x})^{\frac{1}{2}}$. At each iteration $k$, the classical DC algorithm (DCA) \cite{dca_2018} linearizes the second component via its affine minorant $h_k(\boldsymbol{x}) = h(\boldsymbol{x}_k) + \langle \nabla h(\boldsymbol{x}_k), \boldsymbol{x} - \boldsymbol{x}_k \rangle$. Solving the resulting strictly convex subproblem yields
	\begin{equation}
		\begin{aligned}
			\boldsymbol{x}_{k+1} &= \arg \min_{\boldsymbol{x} \in \mathbb{R}^n} \;  \left\{g(\boldsymbol{x}) - h_k(\boldsymbol{x})\right\}\\
			&= \arg \min_{\boldsymbol{x} \in \mathbb{R}^n} \;  \left\{ \|\boldsymbol{x}\|^2 - \frac{\langle \boldsymbol{Ax}_k, \boldsymbol{x} \rangle}{(\boldsymbol{x}_k^T \boldsymbol{A} \boldsymbol{x}_k)^{\frac{1}{2}}} \right\} = \frac{\boldsymbol{Ax}_k}{2(\boldsymbol{x}_k^T \boldsymbol{A} \boldsymbol{x}_k)^{\frac{1}{2}}},
		\end{aligned}
		\label{eq:dca}
	\end{equation}
	coinciding with the equivalent gradient descent step of the power method in \cref{eq:gd_update}.
	
	Within the broader majorization-minimization framework, the DCA update in \cref{eq:dca} corresponds to minimizing the global quadratic surrogate:
	\begin{equation*}
		\boldsymbol{x}_{k+1} = \arg \min_{\boldsymbol{x} \in \mathbb{R}^n} \;  f_k(\boldsymbol{x}),
	\end{equation*}
	where
	\begin{equation}
		\label{eq:pm_sur}
		f_k(\boldsymbol{x}) = f(\boldsymbol{x}_k) + \langle \nabla f(\boldsymbol{x}_k), \boldsymbol{x} - \boldsymbol{x}_k \rangle + \|\boldsymbol{x} - \boldsymbol{x}_k\|^2.
	\end{equation}
	This surrogate constitutes a global majorant tangent at $\boldsymbol{x}_k$. Its Hessian imposes an isotropic curvature, enforcing a uniform spectral upper bound on $f$:
	$$
	\nabla^2 f_k(\boldsymbol{x}_k) = 2 \boldsymbol{I} \succeq \nabla^2 f(\boldsymbol{x}_k).
	$$
	
	While computationally efficient, this isotropic surrogate is inherently loose. This structural limitation motivates the primary inquiry of our work:
	
	\begin{quote}
		\textit{Can we construct a tighter, curvature-aware local surrogate that accelerates convergence, without sacrificing the matrix-free efficiency of the classical power iteration?}
	\end{quote}
	Addressing this trade-off forms the core methodology of the Split-Merge algorithm.
	}
	
	\subsection{Splitting}
	
	To construct a tighter surrogate at the current iterate, we first recall a fundamental decomposition property of PSD matrices.
	
	\begin{lemma}[{\cite[Theorem 7.2.7]{horn2012matrix}}] \label{lem:psd_d}
		A symmetric matrix $\boldsymbol{A} \in \mathbb{R}^{n \times n}$ of rank $r$ is PSD if and only if there exists a full row rank matrix $\boldsymbol{F} \in \mathbb{R}^{r \times n}$ such that
		\begin{equation}
			\boldsymbol{A} = \boldsymbol{F}^T \boldsymbol{F}.
			\label{eq:psd_decomp}
		\end{equation}
	\end{lemma}
	
	\begin{remark}
		The proposed algorithm utilizes the decomposition property of the PSD matrix $\boldsymbol{A}$ only implicitly; an explicit computation of $\boldsymbol{F}$ is not required.
	\end{remark}
	
	For given vectors $\boldsymbol{u}, \boldsymbol{v} \in \mathbb{R}^r$, we define the matrix
	$$
	\boldsymbol{H}_{\boldsymbol{x}} (\boldsymbol{u}, \boldsymbol{v}) = 2 \boldsymbol{I} - \frac{1}{\left(\boldsymbol{x}^T\boldsymbol{A}\boldsymbol{x}\right)^{\frac{1}{2}}} \boldsymbol{F}^T \left(\boldsymbol{u}\boldsymbol{u}^T + \boldsymbol{v} \boldsymbol{v}^T\right)\boldsymbol{F} + \frac{\left(\boldsymbol{Ax}\right) \left(\boldsymbol{Ax}\right)^T}{\left(\boldsymbol{x}^T\boldsymbol{A}\boldsymbol{x}\right)^{\frac{3}{2}}}.
	$$
	
	\begin{proposition} \label{prop:H_succeq}
		Consider a PSD matrix $\boldsymbol{A} \in \mathbb{R}^{n \times n}$ and its factorization $\boldsymbol{A} = \boldsymbol{F}^T \boldsymbol{F}$. For any vectors $\boldsymbol{u}, \boldsymbol{v} \in \mathbb{R}^r$ satisfying $\|\boldsymbol{u}\| \le 1$, $\|\boldsymbol{v}\| \le 1$, and $\boldsymbol{u}^T\boldsymbol{v} = 0$, we have
		\begin{equation}
			\boldsymbol{H}_{\boldsymbol{x}} (\boldsymbol{u}, \boldsymbol{v}) \succeq \nabla^2 f(\boldsymbol{x}), \quad \forall \boldsymbol{x} \in \Theta.
			\label{eq:H_succeq}
		\end{equation}
	\end{proposition}
	
	\begin{proof} 
		The imposed orthogonality and norm constraints directly guarantee the matrix inequality
		$$
		\boldsymbol{I} \succeq \boldsymbol{u}\boldsymbol{u}^T + \boldsymbol{v}\boldsymbol{v}^T.
		$$
		Given the factorized form $\boldsymbol{A} = \boldsymbol{F}^T \boldsymbol{F}$, it follows that
		$$
		\boldsymbol{A} \succeq \boldsymbol{F}^T \left(\boldsymbol{u}\boldsymbol{u}^T+\boldsymbol{v}\boldsymbol{v}^T\right) \boldsymbol{F}.
		$$
		Substituting this bound into the exact Hessian in \cref{eq:hesse_f} yields \cref{eq:H_succeq}.
	\end{proof}
	
	Based on \cref{prop:H_succeq}, we construct a quadratic surrogate of $f$ at $\boldsymbol{x}_k$ as
	\begin{equation*}
		\phi_k (\boldsymbol{x}) = f(\boldsymbol{x}_k) + \left\langle \nabla f(\boldsymbol{x}_k), \boldsymbol{x} - \boldsymbol{x}_k \right\rangle + \frac{1}{2} (\boldsymbol{x} - \boldsymbol{x}_k)^T \boldsymbol{H}_{\boldsymbol{x}_k}(\boldsymbol{u}_k,\boldsymbol{v}_k) (\boldsymbol{x} - \boldsymbol{x}_k).
	\end{equation*}
	Parameterized by $\boldsymbol{u}_k$ and $\boldsymbol{v}_k$, the exact minimization of the surrogate $\phi_k$ directly yields a generalized family of iterative updates:
	\begin{equation}
		\boldsymbol{x}_{k+1} = \arg \min_{\boldsymbol{x} \in \mathbb{R}^n} \;  \phi_k (\boldsymbol{x}).
		\label{eq:general_v}
	\end{equation}
	Specifically, by setting
	\begin{equation} \label{eq:u_sm}
		\boldsymbol{u}_k = \frac{\boldsymbol{F} \boldsymbol{x}_k}{\|\boldsymbol{F} \boldsymbol{x}_k\|},
	\end{equation}
	and defining the scaling factor $\mu_k = 2(\boldsymbol{x}_k^T \boldsymbol{A} \boldsymbol{x}_k)^{\frac{1}{2}}$, the matrix $\boldsymbol{H}_{\boldsymbol{x}_k} (\boldsymbol{u}_k, \boldsymbol{v}_k)$ reduces to a rank-one perturbation of a scaled identity matrix:
	\begin{equation} \label{eq:H_reduced}
		\boldsymbol{H}_{\boldsymbol{x}_k} (\boldsymbol{u}_k, \boldsymbol{v}_k) = 2 \boldsymbol{I} - \frac{2}{\mu_k} \boldsymbol{F}^T \boldsymbol{v}_k \boldsymbol{v}_k^T \boldsymbol{F}.
	\end{equation}
	It is straightforward to verify that the spectrum of this matrix consists of the eigenvalue $2$ with multiplicity at least $n-1$, alongside a single eigenvalue equal to $2\sigma_k$, where we define
	\begin{equation}
		\sigma_k = 1 - {\|\boldsymbol{F}^T \boldsymbol{v}_k\|^2}/{\mu_k}.
		\label{eq:sigma_k}
	\end{equation}
	
	To ensure the positive definiteness of the surrogate Hessian, we restrict our choice of $\boldsymbol{v}_k$ such that $\sigma_k > 0$. Under this condition, an application of the Sherman--Morrison formula \cite{woodbury_1950} to \eqref{eq:H_reduced} yields its inverse as
	\begin{equation}
		\left(\boldsymbol{H}_{\boldsymbol{x}_k} ({\boldsymbol{u}_k,\boldsymbol{v}_k})\right)^{-1} = \frac{1}{2} \left(\boldsymbol{I} + \frac{1}{{\sigma_k \mu_k} } \boldsymbol{F}^T \boldsymbol{v}_k \boldsymbol{v}_k^T \boldsymbol{F}\right).
		\label{eq:H_inv}
	\end{equation}
	
	Substituting the expressions for the gradient \eqref{eq:grad_f} and the inverse Hessian \eqref{eq:H_inv} into the general iteration \eqref{eq:general_v} yields the following update scheme:
	\begin{equation} \label{eq:general_v_expand}
		\begin{aligned}[b]
			\boldsymbol{x}_{k+1} 
			&= \boldsymbol{x}_k - \left(\boldsymbol{H}_{\boldsymbol{x}_k} (\boldsymbol{u}_k,\boldsymbol{v}_k)\right)^{-1} \nabla f(\boldsymbol{x}_k) \\
			&= \boldsymbol{x}_k - \frac{1}{2} \left( \boldsymbol{I} + \frac{1}{\sigma_k \mu_k} \boldsymbol{F}^T \boldsymbol{v}_k \boldsymbol{v}_k^T \boldsymbol{F} \right) \left( 2\boldsymbol{x}_k - \frac{2}{\mu_k} \boldsymbol{A}\boldsymbol{x}_k \right) \\
			&= \frac{1}{\mu_k} \boldsymbol{A}\boldsymbol{x}_k + \left(\frac{\boldsymbol{v}_k^T \boldsymbol{F} \boldsymbol{Ax}_k}{\sigma_k \mu_k^2} - \frac{\boldsymbol{v}_k^T \boldsymbol{u}_k}{2\sigma_k}\right) \boldsymbol{F}^T {\boldsymbol{v}_k} \\
			&= \frac{1}{{\mu_k}} \boldsymbol{Ax}_k + \frac{\boldsymbol{v}_k^T \boldsymbol{F} \boldsymbol{Ax}_k}{{\sigma_k \mu_k^2}} \boldsymbol{F}^T {\boldsymbol{v}_k}.
		\end{aligned}
	\end{equation}
	The final equality seamlessly holds due to the orthogonality condition $\boldsymbol{u}_k^T\boldsymbol{v}_k = 0$.
	
	\begin{remark}[Reduction to the Power Method] \label{rem:reduce_pm}
		In the absence of the auxiliary vector (i.e., $\boldsymbol{v}_k = \boldsymbol{0}$), the update rule \eqref{eq:general_v_expand} immediately reduces to the classical power method. This algebraic consistency motivates our choice of $\boldsymbol{u}_k$ in \eqref{eq:u_sm}.
	\end{remark}
	
	For a valid choice of $\boldsymbol{v}_k$, the resulting surrogate function provides a tighter quadratic upper bound than the standard power method surrogate given in \cref{eq:pm_sur}. We formalize this property as follows.
	\begin{proposition} \label{prop:tigher_pm}
		Let $\boldsymbol{u}_k$ be defined as in \eqref{eq:u_sm}. For any $\boldsymbol{v}_k$ satisfying the conditions of \cref{prop:H_succeq}, it holds that
		$$
		\nabla^2 f_k(\boldsymbol{x}_k) = 2 \boldsymbol{I} \succeq \boldsymbol{H}_{\boldsymbol{x}_k} (\boldsymbol{u}_k,\boldsymbol{v}_k) \succeq \nabla^2 f(\boldsymbol{x}_k).
		$$
	\end{proposition}
	
	\begin{proof}
		The right inequality is given by \cref{prop:H_succeq}. The left inequality is an immediate consequence of \eqref{eq:H_reduced}, as the subtracted term $2\mu_k^{-1} \boldsymbol{F}^T\boldsymbol{v}_k \boldsymbol{v}_k^T \boldsymbol{F}$ is PSD.
	\end{proof}
	
	While the generalized surrogate function yields a tighter quadratic upper bound than the classical power method, its explicit reliance on the decomposition $\boldsymbol{A} = \boldsymbol{F}^T \boldsymbol{F}$ incurs a significant computational bottleneck in large-scale settings. This limitation motivates the following fundamental question:
	
	\begin{quote}
		\textit{Can we construct the auxiliary vector $\boldsymbol{v}_k$ to dictate the accelerated descent trajectory of the surrogate, entirely circumventing the need for explicit matrix factorizations?}
	\end{quote}
	
	\subsection{Merging}
	
	To address the algorithmic bottleneck identified previously, we dynamically select the auxiliary vector $\boldsymbol{v}_k$ to maximize the single-step descent of the surrogate $\phi_k$. This objective motivates the following nested optimization problem:
	\begin{equation}
		\boldsymbol{v}_k = \arg \min_{\boldsymbol{v}\in \Omega_k} \;  \left\{\min_{\boldsymbol{d} \in \mathbb{R}^n} \phi_k(\boldsymbol{x}_k + \boldsymbol{d})\right\},
		\label{eq:optimal_v}
	\end{equation}
	where $\boldsymbol{d} \in \mathbb{R}^n$ denotes the search direction for the variable update (i.e., $\boldsymbol{x}_{k+1} = \boldsymbol{x}_{k}+ \boldsymbol{d}$). The feasible set for the auxiliary vector is defined as
	$$
	\Omega_k = \left\{ \boldsymbol{v} \in \mathbb{R}^r: \boldsymbol{u}_k^T \boldsymbol{v} = 0, \|\boldsymbol{v}\|^2 = {1}/{\rho_k} \right\},
	$$
	where $\rho_k \geq 1$ is a scaling parameter to guarantee $\sigma_k > 0$.
	
	As established in the following theorem, this computationally demanding nested formulation reduces to maximizing a generalized Rayleigh quotient.
	\begin{proposition} \label{prop:gep}
		Determining the optimal solution to problem \eqref{eq:optimal_v} is equivalent to maximizing the following generalized Rayleigh quotient:
		\begin{equation}
			\boldsymbol{v}_k = \arg \max_{\boldsymbol{v} \in \Omega_k} \;  \frac{\boldsymbol{v}^T \boldsymbol{B}_k \boldsymbol{v}}{\boldsymbol{v}^T \boldsymbol{C}_k \boldsymbol{v}},
			\label{eq:gep}
		\end{equation}
		where $\boldsymbol{B}_k = \boldsymbol{F}\boldsymbol{A}\boldsymbol{x}_k (\boldsymbol{F}\boldsymbol{A}\boldsymbol{x}_k)^T$, and $\boldsymbol{C}_k = \rho_k \boldsymbol{I} - 1/\mu_k \boldsymbol{F}\boldsymbol{F}^T$ satisfies $\boldsymbol{v}^T \boldsymbol{C}_k \boldsymbol{v} = \sigma_k > 0$ for any $\boldsymbol{v} \in \Omega_k$.
	\end{proposition}
	
	\begin{proof}
		For any fixed $\boldsymbol{v} \in \Omega_k$, the condition $\sigma_k > 0$ ensures the strict positive definiteness of $\boldsymbol{H}_{\boldsymbol{x}_k} (\boldsymbol{u}_k,\boldsymbol{v})$. Thus, the inner minimization with respect to $\boldsymbol{d}$ yields the unique solution $\boldsymbol{d}_k = - (\boldsymbol{H}_{\boldsymbol{x}_k} (\boldsymbol{u}_k,\boldsymbol{v}))^{-1} \nabla f(\boldsymbol{x}_k)$.
		
		Substituting this optimal step $\boldsymbol{d}_k$ back into \eqref{eq:optimal_v} and omitting the constant terms, the minimization problem is equivalent to maximizing the following objective:
		$$
		\boldsymbol{v}_k = \arg \max_{\boldsymbol{v} \in \Omega_k} \;  \left\{ \nabla f(\boldsymbol{x}_k)^T \left(\boldsymbol{H}_{\boldsymbol{x}_k} (\boldsymbol{u}_k,\boldsymbol{v})\right)^{-1} \nabla f(\boldsymbol{x}_k) \right\}.
		$$
		By incorporating the gradient \eqref{eq:grad_f} and the inverse Hessian \eqref{eq:H_inv}, and applying the orthogonality constraint $\boldsymbol{u}_k^T\boldsymbol{v}=0$, the objective simplifies to
		$$
		\boldsymbol{v}_k = \arg \max_{\boldsymbol{v} \in \Omega_k} \;  \frac{\boldsymbol{v}^T \boldsymbol{B}_k \boldsymbol{v}}{1 - \frac{\boldsymbol{v}^T \boldsymbol{FF}^T \boldsymbol{v}}{\mu_k}}.
		$$
		Finally, utilizing the norm constraint $1 = \rho_k \|\boldsymbol{v}\|^2$ to rewrite the denominator as $\boldsymbol{v}^T \boldsymbol{C}_k \boldsymbol{v}$, we immediately establish the generalized Rayleigh quotient \eqref{eq:gep}.
	\end{proof}
	
	Directly maximizing the generalized Rayleigh quotient in \cref{eq:gep} remains computationally prohibitive in high-dimensional regimes. It intrinsically relies on the explicit matrix factorization $\boldsymbol{A} = \boldsymbol{F}^T\boldsymbol{F}$ and necessitates solving a generalized eigenvalue problem at every iteration. To circumvent this structural bottleneck, we extract a tractable relaxation by bounding the denominator. Since $\boldsymbol{A}$ and $\boldsymbol{F}\boldsymbol{F}^T$ share the same non-zero eigenvalues, the feasible norm constraint $\|\boldsymbol{v}\|^2 = 1/\rho_k$ guarantees that $0 \leq \boldsymbol{v}^T \boldsymbol{F}\boldsymbol{F}^T \boldsymbol{v} \leq \lambda_1 / \rho_k$. Applying these spectral bounds to the objective function, we arrive at the following decoupled relaxation:
	\begin{equation*}
		\left(\boldsymbol{v}^T \boldsymbol{F}\boldsymbol{A}\boldsymbol{x}_k\right)^2 \leq \frac{\boldsymbol{v}^T \boldsymbol{B}_k \boldsymbol{v}}{\boldsymbol{v}^T \boldsymbol{C}_k \boldsymbol{v}} \leq \frac{1}{1 - \frac{\lambda_1}{\rho_k\mu_k}} \left(\boldsymbol{v}^T \boldsymbol{F}\boldsymbol{A}\boldsymbol{x}_k\right)^2.
	\end{equation*}
	
	Following this relaxation, $\boldsymbol{v}_k$ is selected by solving the optimization problem:
	$$
	\boldsymbol{v}_k = \arg \max_{\boldsymbol{v} \in \Omega_k} \; \boldsymbol{v}^T \boldsymbol{F}\boldsymbol{A}\boldsymbol{x}_k \subseteq \arg \max_{\boldsymbol{v} \in \Omega_k} \; \left(\boldsymbol{v}^T \boldsymbol{F}\boldsymbol{A}\boldsymbol{x}_k\right)^2.
	$$
	Recognizing that the feasible set $\Omega_k$ imposes an orthogonality constraint against $\boldsymbol{F}\boldsymbol{x}_k$, this constrained maximization resolves into a pure geometric projection. Its closed-form solution is elegantly captured by a Gram--Schmidt orthogonalization of $\boldsymbol{F}\boldsymbol{A}\boldsymbol{x}_k$ with respect to $\boldsymbol{F}\boldsymbol{x}_k$. Thus, the normalized auxiliary vector is constructed as
	\begin{equation} \label{eq:v_sm}
		\boldsymbol{v}_k = \frac{\boldsymbol{F}\boldsymbol{A}\boldsymbol{x}_k - \alpha_k \boldsymbol{F}\boldsymbol{x}_k}{\|\boldsymbol{F}\boldsymbol{A}\boldsymbol{x}_k - \alpha_k \boldsymbol{F}\boldsymbol{x}_k\|},
	\end{equation}
	where $\alpha_k = (\boldsymbol{x}_k^T \boldsymbol{A}^2 \boldsymbol{x}_k) / (\boldsymbol{x}_k^T \boldsymbol{A} \boldsymbol{x}_k)$ is the projection coefficient. Substituting the scaled vector $\rho_k^{-1/2} \boldsymbol{v}_k$ into \eqref{eq:general_v_expand} yields the following two-term iteration:
	\begin{equation}
		\boldsymbol{x}_{k+1} = \zeta_k \boldsymbol{A}\boldsymbol{x}_k + \omega_k \boldsymbol{A}^2\boldsymbol{x}_k,
		\label{eq:sm_update}
	\end{equation}
	where the combination coefficients are given by
	\begin{equation*}
	\zeta_k = {1}/{\mu_k} - {\alpha_k}/{(\mu_k^2 \sigma_k \rho_k)}, \quad \text{and} \quad \omega_k = {1}/{(\mu_k^2 \sigma_k \rho_k)}.
	\end{equation*}
	
	To simplify the evaluation of $\sigma_k$ in \cref{eq:sigma_k}, we define the residual vector $\boldsymbol{r}_k = \boldsymbol{A}^2\boldsymbol{x}_k - \alpha_k \boldsymbol{A}\boldsymbol{x}_k$ alongside the auxiliary scalar estimator $\gamma_k = {\|\boldsymbol{r}_k\|^2}/{(\boldsymbol{A}\boldsymbol{x}_k)^T \boldsymbol{r}_k}$. The parameter $\sigma_k$ then reduces to the following form:
	$$
	\sigma_k = 1 - {\gamma_k}/{(\rho_k \mu_k)}.
	$$
	
	Crucially, the decomposition factors $\boldsymbol{F}$ and $\boldsymbol{F}^T$ have completely vanished from the final evaluations. The design of our approach captures this profound algebraic interplay: we initially \textit{split} the PSD matrix $\boldsymbol{A}$ to construct a generalized surrogate geometry; subsequently, extracting the auxiliary vector $\boldsymbol{v}$ induces an exact cancellation that intrinsically \textit{merges} these factors back together. This elegant mechanism bypasses explicit matrix factorizations entirely, leading to the matrix-free \textit{Split-Merge} algorithm outlined in \cref{alg:sm_stable}.
	
	{
	\tp
	To establish the theoretical completeness of this algorithmic framework, we formalize the adaptive selection strategy for the scaling factor $\rho_k$ as follows:
	\begin{equation}
		\label{eq:rho_k}
		\rho_k = \max \left\{ 1, \, \frac{\alpha_k}{\mu_k}, \, \frac{\gamma_k}{\mu_k(1-\epsilon_{\sigma})} \right\}.
	\end{equation}
	Here, the parameter $\epsilon_{\sigma} > 0$ (e.g., $10^{-10}$) serves as a safeguard, guaranteeing the positive definiteness condition $\sigma_k \ge \epsilon_{\sigma} > 0$.
	
	The formulation of $\rho_k$ in \cref{eq:rho_k} integrates three spectral estimators: $\alpha_k$, $\gamma_k$, and $\mu_k$. This adaptive mechanism drives $\rho_k \to 1$ to asymptotically tighten the bound. Under transient fluctuations or spectral ill-conditioning, the scaling factor inherently increases to enforce an algebraic safeguard. We formally detail the rigorous convergence to the global minimizer in \cref{sec:convergence}.
	}

	\begin{algorithm}[t]
	\caption{Split-Merge Algorithm}
	\label{alg:sm_stable}
	\begin{algorithmic}[1]
		\tp
		\REQUIRE Symmetric PSD matrix $\boldsymbol{A} \in \mathbb{R}^{n\times n}$, tolerance $\epsilon > 0$, safeguard $\epsilon_{\sigma} > 0$
		\ENSURE Approximation of the dominant eigenvector $\boldsymbol{x}_k$
		\STATE \textbf{Initialize:} Vector $\boldsymbol{x}_0 \in \mathbb{R}^n$ satisfying $\boldsymbol{x}_0^T\boldsymbol{A}\boldsymbol{x}_0 > 0$, and set index $k = 0$
		\WHILE{stopping criterion is not met}
		\STATE $\mu_k = 2 (\boldsymbol{x}_k^T \boldsymbol{A}\boldsymbol{x}_k)^{\frac{1}{2}}$ \COMMENT{Principal magnitude estimator ($\mu_k \to \lambda_1$)}
		\STATE $\alpha_k = \frac{\boldsymbol{x}_k^T \boldsymbol{A}^2 \boldsymbol{x}_k}{\boldsymbol{x}_k^T \boldsymbol{A}\boldsymbol{x}_k}$ \COMMENT{Rayleigh quotient ($\alpha_k \to \lambda_1$)}
		\STATE $\boldsymbol{r}_k = \boldsymbol{A}^2\boldsymbol{x}_k - \alpha_k \boldsymbol{A}\boldsymbol{x}_k$
		\IF{$\|\boldsymbol{r}_k\| \leq \epsilon$}
		\STATE \textbf{break} 
		\ENDIF
		\STATE $\gamma_k = \frac{\|\boldsymbol{r}_k\|^2}{(\boldsymbol{A}\boldsymbol{x}_k)^T \boldsymbol{r}_k}$ \COMMENT{Asymptotic tracker of the sub-dominant spectrum}
		\STATE $\rho_k = \max \left\{ 1, \, \frac{\alpha_k}{\mu_k}, \, \frac{\gamma_k}{\mu_k(1-\epsilon_{\sigma})} \right\}$ \COMMENT{Asymptotic scaling factor ($\rho_k \to 1$)}
		\STATE $\sigma_k = 1 - \frac{\gamma_k}{\rho_k \mu_k}$
		\STATE $\zeta_k = \frac{1}{\mu_k} - \frac{\alpha_k}{\mu_k^2 \sigma_k \rho_k}$
		\STATE $\omega_k = \frac{1}{\mu_k^2 \sigma_k \rho_k}$
		\STATE $\boldsymbol{x}_{k+1} = \zeta_k \boldsymbol{A}\boldsymbol{x}_k + \omega_k \boldsymbol{A}^2\boldsymbol{x}_k$ \COMMENT{Polynomial update step}
		\STATE $k = k + 1$
		\ENDWHILE
	\end{algorithmic}
\end{algorithm}
	
	{
	\tp
	\section{Global Optimality and Acceleration of the Split-Merge Algorithm} \label{sec:convergence}
	
	To establish convergence to the global minimizer for the Split-Merge algorithm and quantify the acceleration induced by the dynamic spectral peeling mechanism, we specify the following initial condition.
	
	\begin{assumption} \label{ass:initialization_alignment}
		Let $\mathcal{Q}_1 = \operatorname{span}\{\boldsymbol{q}_1, \dots, \boldsymbol{q}_{m_1}\}$ denote the dominant eigenspace. The initial iterate $\boldsymbol{x}_0$ satisfies $\mathcal{P}_{\mathcal{Q}_1}(\boldsymbol{x}_0) \neq \boldsymbol{0}$.
	\end{assumption}
	
	Under \cref{ass:initialization_alignment}, without loss of generality, the orthonormal basis of $\mathcal{Q}_1$ is configured such that $\boldsymbol{q}_1$ aligns with this initial projection:
	\begin{equation}
		\boldsymbol{q}_1 = {\mathcal{P}_{\mathcal{Q}_1}(\boldsymbol{x}_0)}/{\|\mathcal{P}_{\mathcal{Q}_1}(\boldsymbol{x}_0)\|}.
		\label{eq:q1_alignment}
	\end{equation}
	
	The Split-Merge update rule \cref{eq:sm_update} can be equivalently formulated as the application of an adaptive matrix polynomial:
	\begin{equation}
		\boldsymbol{x}_{k+1} = P_k(\boldsymbol{A}) \boldsymbol{x}_k,
		\label{eq:poly_operator}
	\end{equation}
	where $P_k(\lambda) = \omega_k \lambda (\lambda - \tau_k)$. In this context, $P_k(\boldsymbol{A})$ acts as a quadratic spectral filter, which is parameterized by its non-trivial root $\tau_k = -\zeta_k / \omega_k$. This root expands to:
	\begin{equation}
		\tau_k = -\rho_k \mu_k + \alpha_k + \gamma_k.
		\label{eq:tau_root}
	\end{equation}
	
	Substituting the adaptive scaling strategy \cref{eq:rho_k} into \cref{eq:tau_root} yields an upper bound on the root $\tau_k$:
	\begin{equation}
		\tau_k \le -\frac{1}{2}\left(\alpha_k + \frac{\gamma_k}{1-\epsilon_{\sigma}}\right) + \alpha_k + \gamma_k < \frac{1}{2}(\alpha_k + \gamma_k).
		\label{eq:tau_bound}
	\end{equation}
	
	\begin{remark}[Spectral Representation of the Estimators] \label{rem:spectral_convexity}
		By the construction of the orthonormal vectors $\boldsymbol{u}_k$ and $\boldsymbol{v}_k$ in \cref{eq:u_sm} and \cref{eq:v_sm}, the scalars $\alpha_k$ and $\gamma_k$ admit the implicit representations $\alpha_k = \boldsymbol{u}_k^T \boldsymbol{F}\boldsymbol{F}^T \boldsymbol{u}_k$ and $\gamma_k = \boldsymbol{v}_k^T \boldsymbol{F}\boldsymbol{F}^T \boldsymbol{v}_k$. Substituting the matrix decomposition $\boldsymbol{F} = \operatorname{diag}(\sqrt{\lambda_1}, \dots, \sqrt{\lambda_r})[\boldsymbol{q}_1, \dots, \boldsymbol{q}_r]^T$, it follows that both $\alpha_k$ and $\gamma_k$ are convex combinations of the non-zero eigenvalues of $\boldsymbol{A}$:
		\begin{equation}
			\alpha_k = \sum_{i=1}^r (\boldsymbol{u}_k)_i^2 \lambda_i, \quad \text{and} \quad \gamma_k = \sum_{i=1}^r (\boldsymbol{v}_k)_i^2 \lambda_i,
			\label{eq:convex_comb}
		\end{equation}
		where the respective convex weights are given by
		\begin{equation}
			(\boldsymbol{u}_k)_i^2 = \frac{\lambda_i(\boldsymbol{q}_i^T\boldsymbol{x}_k)^2}{ \sum_{j=1}^{r} \lambda_j(\boldsymbol{q}_j^T\boldsymbol{x}_k)^2}, \quad \text{and} \quad (\boldsymbol{v}_k)_i^2 = \frac{\lambda_i(\lambda_i - \alpha_k)^2(\boldsymbol{q}_i^T\boldsymbol{x}_k)^2}{\sum_{j=1}^{r} \lambda_j(\lambda_j - \alpha_k)^2(\boldsymbol{q}_j^T\boldsymbol{x}_k)^2}.
			\label{eq:explicit_weights}
		\end{equation}
	\end{remark}
	
	The exact orthogonality $\boldsymbol{u}_k^T \boldsymbol{v}_k = 0$ yields an upper bound on the dynamic root $\tau_k$, thereby driving the attenuation of sub-dominant eigenspaces. To quantify this behavior, we define the attenuation factor:
	\begin{equation}
		\mathcal{A}_k(\lambda) = \left| \frac{\lambda - \tau_k}{\lambda_1 - \tau_k} \right|.
		\label{eq:acceleration_factor}
	\end{equation}
	The following lemma establishes explicit bounds on this factor to characterize the spectral decay.
	
	\begin{lemma}[Sub-dominant Attenuation] \label{lem:acceleration_factor}
		Under \cref{ass:initialization_alignment}, for the sequence $\{\boldsymbol{x}_k\}$ generated by \cref{alg:sm_stable}, the non-trivial root of the spectral filter is bounded as
		\begin{equation}
			\tau_k < \frac{1}{2}(\lambda_1 + \lambda_{m_1+1}).
			\label{eq:tau_strict_bound}
		\end{equation}
		Consequently, the following properties hold:
		
		(i) The iterate $\boldsymbol{x}_k$ maintains a non-zero projection onto the principal eigenvector, i.e., $\boldsymbol{q}_1^T\boldsymbol{x}_{k} \neq 0$ for all $k \ge 0$;
		
		(ii) Within the initial sub-dominant eigenspace, the relative projection ratio decays at least linearly, governed by:
		\begin{equation}
			\left|\frac{\boldsymbol{q}_{i}^T\boldsymbol{x}_{k+1}}{\boldsymbol{q}_1^T\boldsymbol{x}_{k+1}}\right| = \frac{\lambda_{m_1+1}}{\lambda_1} \mathcal{A}_k(\lambda_{m_1+1}) \left|\frac{\boldsymbol{q}_{i}^T\boldsymbol{x}_k}{\boldsymbol{q}_1^T\boldsymbol{x}_k}\right|, \quad \forall i=m_1+1, \dots, m_2,
			\label{eq:spectral_decay}
		\end{equation}
		where the bound $\mathcal{A}_k(\lambda_{m_1+1}) < 1$ ensures faster attenuation relative to the standard power method.
	\end{lemma}
	
	\begin{proof}
		By the polynomial update \cref{eq:poly_operator}, the projections of the iterates onto the eigenvectors $\boldsymbol{q}_i$ evolve as:
		\begin{equation}
			\boldsymbol{q}_i^T\boldsymbol{x}_{k+1} = \omega_k \lambda_i (\lambda_i - \tau_k) \boldsymbol{q}_i^T\boldsymbol{x}_{k}, \quad \forall i = 1, \dots, r, \quad \forall k \ge 0.
			\label{eq:proj_update}
		\end{equation}
		Under \cref{ass:initialization_alignment} and $\boldsymbol{q}_1^T\boldsymbol{x}_{0} > 0$ from \cref{eq:q1_alignment}, the recursion \cref{eq:proj_update} directly enforces:
		\begin{equation}
			\boldsymbol{q}_i^T\boldsymbol{x}_{k} = 0, \quad \forall i \in \{2, \dots, m_1\} \cup \{r+1, \dots, n\}, \quad \forall k \ge 1.
			\label{eq:zero_projections}
		\end{equation}
		Consequently, the associated convex weights \cref{eq:explicit_weights} vanish for these inactive indices:
		\begin{equation}
			(\boldsymbol{u}_k)_i^2 = (\boldsymbol{v}_k)_i^2 = 0, \quad \forall i \in \{2, \dots, m_1\} \cup \{r+1, \dots, n\}, \quad \forall k \ge 1.
			\label{eq:zero_weights}
		\end{equation}
		The orthonormality of the auxiliary vectors $\boldsymbol{u}_k$ and $\boldsymbol{v}_k$ imposes capacity constraints on their spectral components:
		\begin{equation}
			\sum_{i=1}^r \left[ (\boldsymbol{u}_k)_i^2 + (\boldsymbol{v}_k)_i^2 \right] = 2, \quad \text{and} \quad (\boldsymbol{u}_k)_i^2 + (\boldsymbol{v}_k)_i^2 \le 1, \quad \forall i = 1, \dots, r.
			\label{eq:uv_constraints}
		\end{equation}
		Substituting the convex representations \cref{eq:convex_comb} into the upper bound \cref{eq:tau_bound}, and maximizing this spectral combination subject to the zero-weight and capacity constraints \cref{eq:zero_weights,eq:uv_constraints}, yields:
		\begin{equation*}
			\tau_k < \frac{1}{2} \sum_{i=1}^r \left[ (\boldsymbol{u}_k)_i^2 + (\boldsymbol{v}_k)_i^2 \right] \lambda_i \le \frac{1}{2}(\lambda_1 + \lambda_{m_1+1}),
		\end{equation*}
		which immediately establishes \cref{eq:tau_strict_bound}.
		
		Given the strict eigen-gap $\lambda_1 > \lambda_{m_1+1}$, the established bound implies $\tau_k < \lambda_1$. Since $\omega_k > 0$, the principal scaling factor is positive: $\omega_k \lambda_1 (\lambda_1 - \tau_k) > 0$. Combined with the initial condition $\boldsymbol{q}_1^T\boldsymbol{x}_0 \neq 0$, the recursion \cref{eq:proj_update} guarantees the principal projection never vanishes, proving (i).
		
		For (ii), dividing the recursion \cref{eq:proj_update} for $i \in \{m_1+1, \dots, m_2\}$ by the principal recursion directly yields \cref{eq:spectral_decay}. The bound \cref{eq:tau_strict_bound} enforces:
		\begin{equation*}
			\mathcal{A}_k(\lambda_{m_1+1}) = {|\lambda_{m_1+1} - \tau_k|}/{(\lambda_1 - \tau_k)} < 1,
		\end{equation*}
		completing the proof.
	\end{proof}
	
	\begin{remark}[Numerical Well-definedness and Termination] \label{rem:well_definedness}
		\begin{enumerate}
			\item[(i)] \textit{Global Well-definedness}: \cref{lem:acceleration_factor} (i) guarantees $\boldsymbol{x}_k^T \boldsymbol{A} \boldsymbol{x}_k > 0$ for all $k \ge 0$. Consequently, all dynamic parameters ($\alpha_k, \rho_k, \sigma_k, \zeta_k, \omega_k$) in \cref{alg:sm_stable} remain well-defined, precluding any singular behavior in the update rule.
			
			\item[(ii)] \textit{Asymptotic Stability}: As $\boldsymbol{x}_k$ converges to the dominant eigenspace while remaining bounded, the residual vanishes ($\|\boldsymbol{r}_k\| \to 0$). Although the formula for $\gamma_k$ exhibits an indeterminate ratio in the limit, its representation \cref{eq:convex_comb} guarantees its asymptotic confinement within the convex hull of the sub-dominant spectrum. Therefore, the residual norm $\|\boldsymbol{r}_k\|$ serves as a natural stopping criterion, safely terminating the iteration before vanishing terms cause numerical underflow in finite-precision arithmetic.
		\end{enumerate}
	\end{remark}
	
	Building upon the initial sub-dominant attenuation in \cref{lem:acceleration_factor}, the following lemma formalizes the sequential isolation of the residual spectrum via a dynamic peeling mechanism. Once the intermediate sub-dominant components ($i=m_1+1, \dots, t$) achieve sufficient compression, the dynamic root $\tau_k$ shifts deeper into the spectrum. This shift enforces uniform decay across the subsequent eigenspace block $\{i \ge t+1 \mid \lambda_i = \lambda_{t+1}\}$, driving the progressive annihilation of the remaining spectrum.
	
	\begin{lemma}[Cascading Spectral Peeling] \label{lem:spectral_peeling_cascade}
		Under \cref{ass:initialization_alignment}, for the sequence $\{\boldsymbol{x}_k\}$ generated by \cref{alg:sm_stable}, there exists a threshold $\bar{\delta} \in (0, 1]$ such that for any $\delta \in (0, \bar{\delta})$ and any index $t \ge m_2$ exhibiting a strict eigen-gap $\lambda_t > \lambda_{t+1}$, whenever the intermediate eigenspace projections satisfy the attenuation conditions:
		\begin{align}
			\sum_{j=m_1+1}^{t} (\boldsymbol{q}_j^T\boldsymbol{x}_k)^2 &\le \delta^2 (\boldsymbol{q}_1^T\boldsymbol{x}_k)^2, \label{eq:peeling_cond_1} \\
			\sum_{j=m_1+1}^{t} (\boldsymbol{q}_j^T\boldsymbol{x}_k)^2 &\le \delta^2 \sum_{j=t+1}^{r} (\boldsymbol{q}_j^T\boldsymbol{x}_k)^2, \label{eq:peeling_cond_2}
		\end{align}
		the dynamic root of the spectral filter is bounded by:
		\begin{equation}
			\tau_k < \frac{1}{2}(\lambda_1 + \lambda_{t+1}) + \mathcal{C}_{\text{gap}} \delta^2,
			\label{eq:tau_cascading_bound}
		\end{equation}
		where the constant is defined as $\mathcal{C}_{\text{gap}} = \frac{1}{2} (r-m_1)\left(1+\frac{4n^2\lambda_1(\lambda_1-\lambda_r)^2}{\lambda_r(\lambda_1-\lambda_{m_{2}+1})^2}\right)\lambda_{m_1+1}$.
		
		This bounded root intrinsically dictates two properties of the relative projections:
		
		(i) An exact recursive decay within the intermediate block, governed by the contraction $\mathcal{A}_k(\lambda_i) < 1$:
		\begin{equation}
			\left|\frac{\boldsymbol{q}_{i}^T\boldsymbol{x}_{k+1}}{\boldsymbol{q}_1^T\boldsymbol{x}_{k+1}}\right| = \frac{\lambda_{i}}{\lambda_1} \mathcal{A}_k(\lambda_{i})\left|\frac{\boldsymbol{q}_{i}^T\boldsymbol{x}_{k}}{\boldsymbol{q}_1^T\boldsymbol{x}_{k}}\right|, \quad \forall i = m_1+1, \dots, t.
			\label{eq:peeling_decay_exact}
		\end{equation}
		
		(ii) An attenuation bound across the subsequent eigenspace block:
		\begin{equation}
			\left|\frac{\boldsymbol{q}_{i}^T\boldsymbol{x}_{k+1}}{\boldsymbol{q}_1^T\boldsymbol{x}_{k+1}}\right| < \frac{\lambda_{t}}{\lambda_1} \left|\frac{\boldsymbol{q}_{i}^T\boldsymbol{x}_k}{\boldsymbol{q}_1^T\boldsymbol{x}_k}\right|, \quad \forall i \ge t+1 \text{ with } \lambda_i = \lambda_{t+1}.
			\label{eq:peeling_decay_strict}
		\end{equation}
	\end{lemma}
	
	\begin{proof}
		By \cref{eq:explicit_weights} and \cref{eq:peeling_cond_1}, for any index $i \in \{m_1+1, \dots, t\}$, we have:
		\begin{equation}
			(\boldsymbol{u}_k)_i^2 = \frac{\lambda_i}{\lambda_1} \frac{(\boldsymbol{q}_i^T \boldsymbol{x}_k)^2}{(\boldsymbol{q}_1^T \boldsymbol{x}_k)^2} (\boldsymbol{u}_k)_1^2 \le \delta^2.
			\label{eq:bound_u}
		\end{equation}
		Condition \cref{eq:peeling_cond_2} implies the existence of an index $j_0 \in \{t+1, \dots, r\}$ such that 
		$$
		(\boldsymbol{q}_i^T \boldsymbol{x}_k)^2 \le n^2\delta^2 (\boldsymbol{q}_{j_0}^T \boldsymbol{x}_k)^2.
		$$	
		Furthermore, it holds that
		$$
		\max_{j\in\{1, j_0\}} (\lambda_j-\alpha_k)^2 \ge \frac{1}{4}(\lambda_1-\lambda_{j_0})^2 \ge \frac{1}{4}(\lambda_1-\lambda_{m_2+1})^2,
		$$
		where the first inequality follows from bounding the maximum squared distance from below by its value at the midpoint $\alpha_k = (\lambda_1+\lambda_{j_0})/2$. Thus, by \cref{eq:explicit_weights}, selecting $j \in \{1, j_0\}$ yields:
		\begin{equation}
			(\boldsymbol{v}_k)_i^2 = \frac{\lambda_i(\lambda_i-\alpha_k)^2(\boldsymbol{q}_i^T \boldsymbol{x}_k)^2}{\lambda_j(\lambda_j-\alpha_k)^2(\boldsymbol{q}_j^T \boldsymbol{x}_k)^2} (\boldsymbol{v}_k)_j^2 \le \frac{4 n^2 \lambda_1 (\lambda_1-\lambda_r)^2}{\lambda_r(\lambda_1-\lambda_{m_2+1})^2} \delta^2.
			\label{eq:bound_v}
		\end{equation}
		Substituting bounds \cref{eq:bound_u} and \cref{eq:bound_v} into \cref{eq:tau_bound} with \cref{eq:convex_comb}, and applying properties \cref{eq:zero_weights} and \cref{eq:uv_constraints} to maximize the spectral sum, we obtain:
		\begin{equation*}
			\tau_k < \frac{1}{2}(\lambda_1 + \lambda_{t+1}) + \mathcal{C}_{\text{gap}} \delta^2,
		\end{equation*}
		which establishes \cref{eq:tau_cascading_bound}.
		
		The condition $\mathcal{A}_k(\lambda_i) < 1$ for all $i \in \{m_1+1, \dots, t\}$ is equivalent to
		\begin{equation}
			\tau_k < \frac{1}{2}(\lambda_1+\lambda_t).
			\label{eq:tau_req_t}
		\end{equation}
		In view of \cref{eq:tau_cascading_bound}, we define:
		\begin{equation*}
			\bar{\delta}_1(t) = \sqrt{\frac{\lambda_t - \lambda_{t+1}}{2\mathcal{C}_{\text{gap}}}}.
		\end{equation*}
		For any $\delta \in (0, \bar{\delta}_1(t))$, \cref{eq:tau_req_t} holds, which establishes \cref{eq:peeling_decay_exact}.
		
		Next, the bound \cref{eq:peeling_decay_strict} on the subsequent eigenspace block is equivalent to:
		\begin{equation*}
			-\frac{\lambda_{t}}{\lambda_1} < \frac{\lambda_{t+1} (\lambda_{t+1} - \tau_k)}{\lambda_1 (\lambda_1 - \tau_k)} < \frac{\lambda_{t}}{\lambda_1}.
		\end{equation*}
		Since $\lambda_1 - \tau_k > 0$ by \cref{lem:acceleration_factor} and $\lambda_{t+1} < \lambda_t$, the right-hand inequality trivially holds. The left-hand inequality reduces to:
		\begin{equation}
			\tau_k < \frac{\lambda_{t+1}^2 + \lambda_1\lambda_t}{\lambda_t + \lambda_{t+1}}.
			\label{eq:tau_strict_req}
		\end{equation}
		To satisfy \cref{eq:tau_strict_req} via \cref{eq:tau_cascading_bound}, we define the second threshold:
		\begin{equation*}
			\bar{\delta}_2(t) = \sqrt{\frac{(\lambda_1 - \lambda_{t+1})(\lambda_t - \lambda_{t+1})}{2\mathcal{C}_{\text{gap}}(\lambda_t + \lambda_{t+1})}}.
		\end{equation*}
		For any $\delta \in (0, \bar{\delta}_2(t))$, condition \cref{eq:tau_strict_req} holds.
		
		To decouple the bound from $t$, we minimize over all indices with a strict eigen-gap. Since $n$ is finite, this global threshold remains positive:
		\begin{equation*}
			\bar{\delta} = \min_{t \ge m_2, \, \lambda_t > \lambda_{t+1}} \min\left\{1, \bar{\delta}_1(t), \bar{\delta}_2(t)\right\} > 0.
		\end{equation*}
		Consequently, any $\delta \in (0, \bar{\delta})$ ensures \cref{eq:peeling_decay_exact} and \cref{eq:peeling_decay_strict} hold uniformly for all valid $t$, completing the proof.
	\end{proof}
	
	\begin{remark}[Acceleration via Dynamic Peeling] \label{rem:peeling_acceleration}
		\cref{lem:spectral_peeling_cascade} details the acceleration mechanism of the framework. While the standard power method is restricted by the static decay ratio $\lambda_i/\lambda_1$, \cref{eq:peeling_decay_exact} demonstrates that dynamic spectral peeling breaks this linear bottleneck. As the iteration proceeds, the dynamic root $\tau_k$ shifts into the interior spectrum. Whenever $\tau_k$ approximates a sub-dominant eigenvalue $\lambda_i$, the corresponding factor $\mathcal{A}_k(\lambda_i) = |\lambda_i - \tau_k| / (\lambda_1 - \tau_k)$ approaches zero. This yields a targeted suppression of eigenspaces clustered near $\tau_k$, thereby accelerating convergence.
	\end{remark}
	
	The ensuing lemma establishes that the sequential compression forces the non-dominant spectrum to vanish below any arbitrary tolerance in finitely many iterations.
	
	\begin{lemma}[Global Sub-dominant Suppression] \label{lem:global_suppression}
		Under \cref{ass:initialization_alignment}, for the sequence $\{\boldsymbol{x}_k\}$ generated by \cref{alg:sm_stable}, for any tolerance $\epsilon \in (0, 1)$, there exists a finite iteration index $K \ge 0$ such that:
		\begin{equation}
			\left|\frac{\boldsymbol{q}_i^T\boldsymbol{x}_K}{\boldsymbol{q}_1^T\boldsymbol{x}_K}\right| < \epsilon, \quad \forall i \in \{m_1+1, \dots, r\}.
			\label{eq:global_epsilon_bound}
		\end{equation}
	\end{lemma}
	
	\begin{proof}
		Fix $\delta \in (0, \bar{\delta})$ from \cref{lem:spectral_peeling_cascade}, and define $\bar{\epsilon} = \delta \epsilon / \sqrt{r}$.
		By \cref{lem:acceleration_factor}, the relative projections for the initial sub-dominant block decay at least linearly. Thus, there exists an index $K_1 \ge 0$ satisfying:
		\begin{equation*}
			\left|\frac{\boldsymbol{q}_i^T\boldsymbol{x}_k}{\boldsymbol{q}_1^T\boldsymbol{x}_k}\right| \le \bar{\epsilon}, \quad \forall i \in \{m_1+1, \dots, m_2\}, \quad \forall k \ge K_1.
		\end{equation*}
		Proceeding by induction, assume that for some block boundary index $t \ge m_2$ and iteration $K_t \ge 0$:
		\begin{equation}
			\left|\frac{\boldsymbol{q}_i^T\boldsymbol{x}_k}{\boldsymbol{q}_1^T\boldsymbol{x}_k}\right| \le \bar{\epsilon}, \quad \forall i \in \{m_1+1, \dots, t\}, \quad \forall k \ge K_t.
			\label{eq:induction_hypo}
		\end{equation}
		By the definition of $\bar{\epsilon}$ and the hypothesis \cref{eq:induction_hypo}, any $k \ge K_t$ satisfies:
		\begin{equation*}
			\frac{\sum_{j=m_1+1}^t (\boldsymbol{q}_j^T\boldsymbol{x}_k)^2}{(\boldsymbol{q}_1^T\boldsymbol{x}_k)^2} \le (t-m_1)\bar{\epsilon}^2 < r \bar{\epsilon}^2 = \delta^2 \epsilon^2 < \delta^2.
		\end{equation*}
		Consequently, the first peeling condition \cref{eq:peeling_cond_1} holds for all $k \ge K_t$. The second condition \cref{eq:peeling_cond_2} bifurcates into two cases:
		
		\textit{Case (i):} Suppose \cref{eq:peeling_cond_2} fails at some iteration $K \ge K_t$. This implies:
		\begin{equation*}
			\sum_{j=t+1}^r (\boldsymbol{q}_j^T\boldsymbol{x}_K)^2 < \frac{1}{\delta^2} \sum_{j=m_1+1}^t (\boldsymbol{q}_j^T\boldsymbol{x}_K)^2.
		\end{equation*}
		Incorporating the established bound, for any $i \in \{t+1, \dots, r\}$, it directly follows that:
		\begin{equation*}
			\left|\frac{\boldsymbol{q}_i^T\boldsymbol{x}_K}{\boldsymbol{q}_1^T\boldsymbol{x}_K}\right|^2 < \frac{1}{\delta^2} \frac{\sum_{j=m_1+1}^t (\boldsymbol{q}_j^T\boldsymbol{x}_K)^2}{(\boldsymbol{q}_1^T\boldsymbol{x}_K)^2} \le \frac{t-m_1}{\delta^2} \bar{\epsilon}^2 < \frac{r}{\delta^2} \bar{\epsilon}^2 = \epsilon^2.
		\end{equation*}
		Since $\bar{\epsilon} < \epsilon$, the hypothesis \cref{eq:induction_hypo} and the above inequality jointly guarantee that \cref{eq:global_epsilon_bound} strictly holds for all $i \in \{m_1+1, \dots, r\}$ at iteration $K$, abruptly terminating the induction.
		
		\textit{Case (ii):} Suppose \cref{eq:peeling_cond_2} holds for all $k \ge K_t$. By \cref{lem:spectral_peeling_cascade}, the decay bounds \cref{eq:peeling_decay_exact} and \cref{eq:peeling_decay_strict} activate. Let $t' \ge t+1$ denote the terminal index of the next eigenspace block (the largest index with $\lambda_{t'} = \lambda_{t+1}$). These bounds jointly enforce a contraction uniformly across all indices up to this block:
		\begin{equation*}
			\left|\frac{\boldsymbol{q}_i^T\boldsymbol{x}_{k+1}}{\boldsymbol{q}_1^T\boldsymbol{x}_{k+1}}\right| < \frac{\lambda_{m_1+1}}{\lambda_1} \left|\frac{\boldsymbol{q}_i^T\boldsymbol{x}_k}{\boldsymbol{q}_1^T\boldsymbol{x}_k}\right|, \quad \forall i \in \{m_1+1, \dots, t'\}, \quad \forall k \ge K_t.
		\end{equation*}
		Thus, there exists an iteration $K_{t'} \ge K_t$ such that this ratio drops below $\bar{\epsilon}$ for all $k \ge K_{t'}$, extending the induction hypothesis to the block boundary $t'$.
		
		Because the rank $r$ is finite, this progressive block-wise peeling must terminate. Thus, in either case, \cref{eq:global_epsilon_bound} is achieved in finite iterations, completing the proof.
	\end{proof}
	
	The preceding sequence of sub-dominant suppressions ultimately culminates in the convergence to a global minimizer of $f$.
	
	\begin{theorem}[Convergence to the Global Minimizer] \label{thm:global_convergence}
		Under \cref{ass:initialization_alignment}, the sequence $\{\boldsymbol{x}_k\}$ generated by \cref{alg:sm_stable} converges to a global minimizer of the objective $f$, satisfying:
		\begin{equation} \label{eq:asymptotic_limit}
			\lim_{k \to \infty} \boldsymbol{x}_k = \sqrt{\lambda_1} \boldsymbol{q}_1 / 2.
		\end{equation}
	\end{theorem}
	
	\begin{proof}
		Define the normalized iterate $\boldsymbol{p}_k = \boldsymbol{x}_k / \|\boldsymbol{x}_k\|$. We first establish the directional convergence:
		\begin{equation}
			\lim_{k \to \infty} \boldsymbol{p}_k = \boldsymbol{q}_1.
			\label{eq:directional_limit}
		\end{equation}
		By \cref{eq:zero_projections}, establishing \cref{eq:directional_limit} reduces to proving the asymptotic vanishing of the remaining sub-dominant projections:
		\begin{equation} 
			\lim_{k \to \infty} \left|\frac{\boldsymbol{q}_i^T\boldsymbol{x}_{k}}{\boldsymbol{q}_1^T\boldsymbol{x}_{k}}\right| = 0, \quad \forall i = m_1+1, \dots, r.
			\label{eq:subdominant_limit}
		\end{equation}
		By \cref{lem:global_suppression}, for any tolerance $\epsilon > 0$, there exists an index $K(\epsilon)$ such that \cref{eq:global_epsilon_bound} holds. Fix $\delta \in (0, \bar{\delta})$ from \cref{lem:spectral_peeling_cascade}, and define the constant $L = (\lambda_1 + \lambda_{m_1+1}) / (\lambda_1 - \lambda_{m_1+1}) > 1$. Construct the auxiliary sequence:
		\begin{equation}
			c_{m_1+1} = \dots = c_{m_2} = 1, \quad c_j = \frac{L}{\delta} \sum_{i=m_1+1}^{j-1} c_i, \quad \forall j = m_2+1, \dots, r.
			\label{eq:aux_c}
		\end{equation}
		By construction, these constants satisfy:
		\begin{equation}
			c_i \ge 1, \quad \forall i = m_1+1, \dots, r.
			\label{eq:c_lower_bound}
		\end{equation}
		Setting $\epsilon = \delta / \sum_{i=m_1+1}^r c_i$, we claim the relative projections are uniformly bounded:
		\begin{equation}
			\left|\frac{\boldsymbol{q}_i^T\boldsymbol{x}_k}{\boldsymbol{q}_1^T\boldsymbol{x}_k}\right| \le c_i \epsilon, \quad \forall i = m_1+1, \dots, r, \quad \forall k \ge K(\epsilon).
			\label{eq:claim_c_bound}
		\end{equation}
		Since the constants $c_i$ are independent of $k$, the bound \cref{eq:claim_c_bound} directly implies \cref{eq:subdominant_limit}, thereby establishing \cref{eq:directional_limit}.
		
		We proceed by contradiction. Suppose \cref{eq:claim_c_bound} fails. By \cref{lem:acceleration_factor}, the initial sub-dominant block decays at least linearly, satisfying the bound for sufficiently large $k$. Thus, the violation must occur at some index $i_0 \in \{m_2+1, \dots, r\}$. Specifically, there exists $k > K(\epsilon)$ such that $\left|\boldsymbol{q}_{i_0}^T\boldsymbol{x}_k / \boldsymbol{q}_1^T\boldsymbol{x}_k\right| > c_{i_0}\epsilon$, while the bound holds for all preceding indices:
		\begin{equation}
			\left|\frac{\boldsymbol{q}_i^T\boldsymbol{x}_k}{\boldsymbol{q}_1^T\boldsymbol{x}_k}\right| \le c_i \epsilon, \quad \forall i \in \{m_1+1, \dots, i_0-1\}, \quad \forall k \ge K(\epsilon).
			\label{eq:hypo_preceding}
		\end{equation}
		By \cref{eq:tau_strict_bound}, $\mathcal{A}_{k}(\lambda_{i}) < L$ for all $i \ge m_1+1$. Inverting the recursion \cref{eq:proj_update} for $i_0$ yields:
		\begin{equation}
			\left|\frac{\boldsymbol{q}_{i_0}^T\boldsymbol{x}_{k-1}}{\boldsymbol{q}_1^T\boldsymbol{x}_{k-1}}\right| = \frac{1}{\mathcal{A}_{k}(\lambda_{i_0})} \frac{\lambda_1}{\lambda_{i_0}} \left|\frac{\boldsymbol{q}_{i_0}^T\boldsymbol{x}_k}{\boldsymbol{q}_1^T\boldsymbol{x}_k}\right| > \frac{c_{i_0}\epsilon}{L}.
			\label{eq:amplification_bound}
		\end{equation}
		By \cref{eq:aux_c} and \cref{eq:hypo_preceding}, the spectral peeling conditions \cref{eq:peeling_cond_1,eq:peeling_cond_2} are satisfied at iteration $k-1$:
		\begin{align*}
			\frac{\sum_{i=m_1+1}^{i_0-1} (\boldsymbol{q}_i^T\boldsymbol{x}_{k-1})^2}{(\boldsymbol{q}_1^T\boldsymbol{x}_{k-1})^2} &\le \left( \sum_{i=m_1+1}^{i_0-1} c_i \epsilon \right)^2 \le \delta^2, \\
			\frac{\sum_{i=m_1+1}^{i_0-1} (\boldsymbol{q}_i^T\boldsymbol{x}_{k-1})^2}{(\boldsymbol{q}_{i_0}^T\boldsymbol{x}_{k-1})^2} &< \frac{\left(\sum_{i=m_1+1}^{i_0-1} c_i \epsilon\right)^2}{(c_{i_0} \epsilon / L)^2} = \delta^2.
		\end{align*}
		By \cref{lem:spectral_peeling_cascade}, this induces the decay in \cref{eq:peeling_decay_strict}, enforcing:
		\begin{equation*}
			\left|\frac{\boldsymbol{q}_{i_0}^T\boldsymbol{x}_{k-1}}{\boldsymbol{q}_1^T\boldsymbol{x}_{k-1}}\right| > \frac{\lambda_1}{\lambda_{i_0-1}} \left|\frac{\boldsymbol{q}_{i_0}^T\boldsymbol{x}_k}{\boldsymbol{q}_1^T\boldsymbol{x}_k}\right| > c_{i_0} \epsilon.
		\end{equation*}
		Applying this cascading argument recursively back to $K(\epsilon)$ yields the chain:
		\begin{equation*}
			c_{i_0} \epsilon < \left|\frac{\boldsymbol{q}_{i_0}^T\boldsymbol{x}_k}{\boldsymbol{q}_1^T\boldsymbol{x}_k}\right| < \left|\frac{\boldsymbol{q}_{i_0}^T\boldsymbol{x}_{k-1}}{\boldsymbol{q}_1^T\boldsymbol{x}_{k-1}}\right| < \dots < \left|\frac{\boldsymbol{q}_{i_0}^T\boldsymbol{x}_{K(\epsilon)}}{\boldsymbol{q}_1^T\boldsymbol{x}_{K(\epsilon)}}\right| < \epsilon,
		\end{equation*}
		which enforces $c_{i_0} < 1$, contradicting \cref{eq:c_lower_bound}. Thus, \cref{eq:claim_c_bound} holds.
		
		Next, we determine the limit of the magnitude $t_k = \|\boldsymbol{x}_k\|$. Substituting $\boldsymbol{x}_k = t_k \boldsymbol{p}_k$ into the update \cref{eq:sm_update} gives:
		\begin{equation*}
			\begin{aligned}
				\boldsymbol{x}_{k+1} &= \left[\frac{1}{2(\boldsymbol{x}_k^T \boldsymbol{A} \boldsymbol{x}_k)^{1/2}} - \frac{\boldsymbol{x}_k^T\boldsymbol{A}^2\boldsymbol{x}_k}{4\sigma_k\rho_k(\boldsymbol{x}_k^T \boldsymbol{A} \boldsymbol{x}_k)^2} \right]\boldsymbol{A}\boldsymbol{x}_k + \frac{1}{4\sigma_k\rho_k(\boldsymbol{x}_k^T \boldsymbol{A} \boldsymbol{x}_k)} \boldsymbol{A}^2\boldsymbol{x}_k \\
				&= \frac{\boldsymbol{A}\boldsymbol{p}_k}{2(\boldsymbol{p}_k^T \boldsymbol{A} \boldsymbol{p}_k)^{1/2}} - \frac{(\boldsymbol{p}_k^T\boldsymbol{A}^2\boldsymbol{p}_k) \boldsymbol{A}\boldsymbol{p}_k }{4\sigma_k\rho_k t_k(\boldsymbol{p}_k^T \boldsymbol{A} \boldsymbol{p}_k)^2} + \frac{\boldsymbol{A}^2\boldsymbol{p}_k}{4\sigma_k\rho_k t_k(\boldsymbol{p}_k^T \boldsymbol{A} \boldsymbol{p}_k)}.
			\end{aligned}
		\end{equation*}
		Left-multiplying by $\boldsymbol{p}_k^T$ exactly cancels the latter two terms, yielding:
		\begin{equation*}
			t_{k+1}\boldsymbol{p}_k^T\boldsymbol{p}_{k+1} = {(\boldsymbol{p}_k^T\boldsymbol{A}\boldsymbol{p}_k)^{1/2}}/{2}. 
		\end{equation*}
		Taking the limit $k \to \infty$, the directional convergence \cref{eq:directional_limit} ensures:
		\begin{equation*}
			\lim_{k \to \infty} \boldsymbol{p}_k^T\boldsymbol{p}_{k+1} = 1, \quad \text{and} \quad \lim_{k \to \infty} \boldsymbol{p}_k^T\boldsymbol{A}\boldsymbol{p}_k = \lambda_1.
		\end{equation*}
		Consequently, the sequence of magnitudes converges:
		\begin{equation} \label{eq:magnitude_limit}
			\lim_{k \to \infty} t_k = {\sqrt{\lambda_1}}/{2}.
		\end{equation}
		Combining \cref{eq:directional_limit} and \cref{eq:magnitude_limit} establishes \cref{eq:asymptotic_limit}. By \cref{thm:optimality}, this limit point is a global minimizer of $f$, completing the proof.
	\end{proof}
	
	\begin{remark}[Adaptive Spectral Extraction] \label{rem:adaptive_extraction}
		As $\boldsymbol{x}_k \to \sqrt{\lambda_1} \boldsymbol{q}_1 / 2$, the principal estimators satisfy $\mu_k \to \lambda_1$ and $\alpha_k \to \lambda_1$. Concurrently, $\gamma_k$ asymptotically reduces to a convex combination of the sub-dominant eigenvalues: $\gamma_k - \sum_{i=m_1+1}^r (\boldsymbol{v}_k)_i^2 \lambda_i \to 0$. Under the scaling formulation \cref{eq:rho_k}, these limits ensure $\rho_k \to 1$. Consequently, the dynamic root \cref{eq:tau_root} asymptotically aligns with $\gamma_k$, yielding $\tau_k - \gamma_k \to 0$. Unlike classical accelerated schemes \cite{ppm, power_m} requiring explicit prior knowledge of the sub-dominant eigenvalues, the Split-Merge algorithm extracts this sub-dominant information directly from the iteration sequence.
	\end{remark}
	
	\begin{figure}[tbp]
		\centering
		\includegraphics[width=0.7\linewidth]{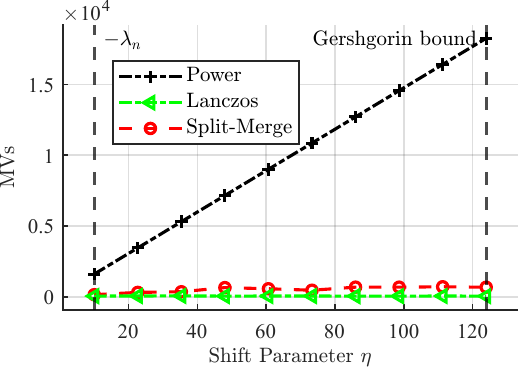}
		\caption{\tp\textbf{Impact of diagonal shifts on computational cost.} The matrix-vector products are evaluated against the shift parameter $\eta$ applied to $\boldsymbol{A} + \eta \boldsymbol{I}$, where $\eta$ ranges from the exact $-\lambda_n$ to the conservative Gershgorin bound.}
		\label{fig:shift_eta_impact}
	\end{figure}
	
	\begin{remark}[Insensitivity to Diagonal Shifts] \label{rem:shift_insensitivity}
		For indefinite matrices, applying a diagonal shift $\boldsymbol{A} + \eta \boldsymbol{I}$ is standard practice. Under this shift, the static decay ratio of the power method $(\lambda_i + \eta)/(\lambda_1 + \eta)$ deteriorates toward $1$ as $\eta \to \infty$. The Split-Merge algorithm prevents this deterioration. Since the dynamic root $\tau_k$ asymptotically aligns with a convex combination of the sub-dominant spectrum, it translates uniformly with the shifted eigenvalues. Consequently, the parameter $\eta$ cancels within the asymptotic attenuation factor \cref{eq:acceleration_factor}. This property ensures a convergence rate that is highly insensitive to the shift parameter $\eta$, thereby preserving the targeted subspace suppression. This robustness is empirically corroborated in \cref{subsec:synthetic}, with the consistent trajectories visualized in \cref{fig:shift_eta_impact}.
	\end{remark}
	}
	
	\section{Experiments} \label{sec:experiment}
	
	In this section, we evaluate the empirical performance of the proposed Split-Merge algorithm\footnote{\tp Available at \url{https://github.com/xzliu-opt/SplitMerge}} across a variety of synthetic and real-world datasets. To demonstrate its numerical efficacy and scalability, we benchmark our approach against a comprehensive suite of established eigensolvers, categorized as follows:
	\begin{itemize}
		\item \textbf{Power iteration techniques:} The standard power method \cite{pm_1929} and its accelerated momentum variant, Power+M\footnote{Available at \url{https://github.com/git-xp/Accelerated_PCA}} \cite{power_m}.
		\item \textbf{Subspace methods:} The Lanczos method \cite{lanczos}, {\tp the LOBPCG method}\footnote{\tp Available at \url{https://github.com/lobpcg/blopex}} \cite{LOBPCG}, and the JD method\footnote{Available at \url{https://webspace.science.uu.nl/~sleij101/index.html}} \cite{jd_2000}.
		\item \textbf{Difference-based approaches:} Newton's method \cite{mongeau2004computing} and the L-BFGS method \cite{shi2016limited}.
		\item {\tp \textbf{Optimized software packages:} MATLAB's native \texttt{eigs} solver. Driven by the ARPACK implementation of the implicitly restarted Lanczos method and executed through compiled C/Fortran routines, this solver serves as a robust, industry-standard benchmark for practical eigenvalue computations.}
	\end{itemize}
	
	All baseline methods are evaluated using their recommended default parameters. {\tp For the subspace methods, specifically LOBPCG and JD, we employ the identity matrix as the default preconditioner.} For the Power+M algorithm, we establish an idealized baseline, denoted as {Power+M (ideal)}, by adopting the theoretical optimal momentum parameter $\beta^* = \lambda_2^2 / 4$. Because this exact configuration requires unobservable spectral knowledge that is rarely available in practice, we also evaluate a perturbed variant with $\beta = 0.9 \beta^*$, denoted as {Power+M (near-ideal)}, to reflect a more realistic scenario. Collectively, these comprehensive comparisons are designed to highlight the superior efficiency and practical advantages of our proposed algorithm.

	For all algorithms, the initial vector $\boldsymbol{x}_0$ is generated using MATLAB's \texttt{randn} function and remains consistent across methods.
	All results are averaged over 500 random trials.
	The stopping criterion is defined as $\sin(\theta_k) \leq \epsilon$, where $\epsilon = 10^{-5}$ and $\theta_k$ is the angle between the current iterate $\boldsymbol{x}_k$  and the dominant eigenvector $\boldsymbol{q}_1$.
	Alternatively, the algorithm terminates if the number of iterations exceeds 20,000.
	
	{\tp Algorithmic performance is evaluated primarily in terms of computation time and the number of matrix-vector products (MVs). However, MV counts are excluded for Newton's method (whose complexity is dominated by linear system solves) and the black-box \texttt{eigs} solver, for which only computation times are reported.}
	
	All experiments were performed in MATLAB R2021b on a Windows system equipped with an Intel Core i7-12700H processor (2.30 GHz) and 16 GB of RAM.
	
	\begin{table}[htbp]
	\centering
	\caption{\tp \textbf{Comparison of average computation time (sec.) and the number of matrix-vector products across various matrix dimensions ($n$) and eigen-gaps ($\Delta$)}. Bold and underlined values indicate the \textbf{best} and \underline{second-best} results, respectively. The proposed method is highlighted as \colorbox{smblue}{Split-Merge}. Additionally, the \textit{speed-up} factor of the Split-Merge algorithm relative to the power method is reported.}
	\label{tab:tabsota}
	\renewcommand{\arraystretch}{1.12}
	\setlength{\tabcolsep}{4pt}
	\footnotesize
	\begin{tabular}{l|cc|cc|cc}
		\toprule
		\multirow{2}{*}{\textbf{Method}} & \textbf{Time} & \textbf{MVs} & \textbf{Time} & \textbf{MVs} & \textbf{Time} & \textbf{MVs} \\
		\cmidrule(lr){2-3}\cmidrule(lr){4-5}\cmidrule(lr){6-7}
		& \multicolumn{2}{c}{\textbf{$\Delta = 10^{-1}$}} & \multicolumn{2}{c}{\textbf{$\Delta = 10^{-2}$}} & \multicolumn{2}{c}{\textbf{$\Delta = 10^{-3}$}} \\
		\midrule
		\multicolumn{7}{c}{\textbf{$n = 1024$}} \\
		\midrule
		\textbf{Power} & $3.80\times10^{-3}$ & 93 & $6.77\times10^{-3}$ & 181 & $2.83\times10^{-2}$ & 777 \\
		\textbf{Power$+$M (ideal)} & $1.37\times10^{-3}$ & 30 & $1.87\times10^{-3}$ & 42 & $3.91\times10^{-3}$ & 91 \\
		\textbf{Power$+$M (near-ideal)} & $1.80\times10^{-3}$ & 40 & $2.92\times10^{-3}$ & 69 & $1.08\times10^{-2}$ & 258 \\
		\textbf{Lanczos} & \underline{$1.10\times10^{-3}$} & \textbf{12} & \underline{$1.47\times10^{-3}$} & \textbf{17} & \underline{$3.77\times10^{-3}$} & \textbf{33} \\
		\textbf{LOBPCG} & $6.85\times10^{-3}$ & \underline{14} & $8.44\times10^{-3}$ & \underline{20} & $1.56\times10^{-2}$ & \underline{39} \\
		\textbf{JD} & $1.17\times10^{-2}$ & 39 & $1.40\times10^{-2}$ & 58 & $2.01\times10^{-2}$ & 91 \\
		\textbf{Newton} & $1.68\times10^{-2}$ & -- & $1.61\times10^{-2}$ & -- & $1.60\times10^{-2}$ & -- \\
		\textbf{L-BFGS} & $3.46\times10^{-3}$ & 71 & $4.20\times10^{-3}$ & 89 & $7.69\times10^{-3}$ & 154 \\
		\rowcolor{smblue}\textbf{Split-Merge} & $\mathbf{8.28\times10^{-4}}$ & 17 & $\mathbf{1.03\times10^{-3}}$ & 25 & $\mathbf{3.62\times10^{-3}}$ & 88 \\
		\midrule
		\texttt{eigs} & $2.19\times10^{-3}$ & -- & $2.20\times10^{-3}$ & -- & $4.12\times10^{-3}$ & -- \\
		\midrule
		\textit{Speed-up} & \multicolumn{2}{c}{\textit{4.59$\times$}} & \multicolumn{2}{c}{\textit{6.58$\times$}} & \multicolumn{2}{c}{\textit{7.81$\times$}} \\
		\midrule[1.2pt]
		\addlinespace[2pt]
		\multicolumn{7}{c}{\textbf{$n = 2048$}} \\
		\midrule
		\textbf{Power} & $2.95\times10^{-2}$ & 100 & $5.70\times10^{-2}$ & 194 & $2.41\times10^{-1}$ & 823 \\
		\textbf{Power$+$M (ideal)} & $9.81\times10^{-3}$ & 31 & $1.40\times10^{-2}$ & 45 & $2.93\times10^{-2}$ & 93 \\
		\textbf{Power$+$M (near-ideal)} & $1.35\times10^{-2}$ & 43 & $2.29\times10^{-2}$ & 74 & $8.54\times10^{-2}$ & 273 \\
		\textbf{Lanczos} & $\mathbf{4.92\times10^{-3}}$ & \textbf{13} & $\mathbf{6.76\times10^{-3}}$ & \textbf{18} & $\mathbf{1.56\times10^{-2}}$ & \textbf{36} \\
		\textbf{LOBPCG} & $1.50\times10^{-2}$ & \underline{15} & $1.93\times10^{-2}$ & \underline{21} & $3.68\times10^{-2}$ & \underline{42} \\
		\textbf{JD} & $2.97\times10^{-2}$ & 41 & $4.03\times10^{-2}$ & 61 & $6.05\times10^{-2}$ & 94 \\
		\textbf{Newton} & $9.02\times10^{-2}$ & -- & $9.01\times10^{-2}$ & -- & $9.07\times10^{-2}$ & -- \\
		\textbf{L-BFGS} & $2.28\times10^{-2}$ & 74 & $3.01\times10^{-2}$ & 94 & $5.54\times10^{-2}$ & 162 \\
		\rowcolor{smblue}\textbf{Split-Merge} & \underline{$5.38\times10^{-3}$} & 18 & \underline{$7.68\times10^{-3}$} & 26 & \underline{$2.71\times10^{-2}$} & 91 \\
		\midrule
		\texttt{eigs} & $8.30\times10^{-3}$ & -- & $9.11\times10^{-3}$ & -- & $1.78\times10^{-2}$ & -- \\
		\midrule
		\textit{Speed-up} & \multicolumn{2}{c}{\textit{5.48$\times$}} & \multicolumn{2}{c}{\textit{7.43$\times$}} & \multicolumn{2}{c}{\textit{8.89$\times$}} \\
		\midrule[1.2pt]
		\addlinespace[2pt]
		\multicolumn{7}{c}{\textbf{$n = 4096$}} \\
		\midrule
		\textbf{Power} & $2.90\times10^{-1}$ & 106 & $5.59\times10^{-1}$ & 207 & $2.36\times10^{0}$ & 873 \\
		\textbf{Power$+$M (ideal)} & $9.15\times10^{-2}$ & 33 & $1.29\times10^{-1}$ & 47 & $2.67\times10^{-1}$ & 97 \\
		\textbf{Power$+$M (near-ideal)} & $1.26\times10^{-1}$ & 46 & $2.16\times10^{-1}$ & 79 & $7.92\times10^{-1}$ & 290 \\
		\textbf{Lanczos} & $\mathbf{3.84\times10^{-2}}$ & \textbf{13} & $\mathbf{5.35\times10^{-2}}$ & \textbf{19} & $\mathbf{1.10\times10^{-1}}$ & \textbf{37} \\
		\textbf{LOBPCG} & $6.89\times10^{-2}$ & \underline{16} & $9.17\times10^{-2}$ & \underline{22} & \underline{$1.75\times10^{-1}$} & \underline{45} \\
		\textbf{JD} & $1.54\times10^{-1}$ & 43 & $2.31\times10^{-1}$ & 68 & $3.24\times10^{-1}$ & 96 \\
		\textbf{Newton} & $3.87\times10^{-1}$ & -- & $3.87\times10^{-1}$ & -- & $3.88\times10^{-1}$ & -- \\
		\textbf{L-BFGS} & $2.02\times10^{-1}$ & 76 & $2.66\times10^{-1}$ & 99 & $4.71\times10^{-1}$ & 168 \\
		\rowcolor{smblue}\textbf{Split-Merge} & \underline{$5.12\times10^{-2}$} & 19 & \underline{$7.25\times10^{-2}$} & 27 & $2.58\times10^{-1}$ & 96 \\
		\midrule
		\texttt{eigs} & $5.87\times10^{-2}$ & -- & $6.54\times10^{-2}$ & -- & $1.27\times10^{-1}$ & -- \\
		\midrule
		\textit{Speed-up} & \multicolumn{2}{c}{\textit{5.67$\times$}} & \multicolumn{2}{c}{\textit{7.72$\times$}} & \multicolumn{2}{c}{\textit{9.16$\times$}} \\
		\bottomrule
	\end{tabular}
\end{table}
	
	\subsection{Synthetic Datasets} \label{subsec:synthetic}
	
	{
	\tp
	We first evaluate the methods on synthetic PSD matrices to analyze performance across varying dimensions ($n$) and eigen-gap values ($\Delta$). The test matrices are generated as $\boldsymbol{A} = \boldsymbol{Q} \operatorname{diag}(\lambda_1, \dots, \lambda_n) \boldsymbol{Q}^T \in \mathbb{R}^{n \times n}$, with the orthogonal matrix $\boldsymbol{Q}$ derived from a standard Gaussian QR factorization. The eigen-gap is controlled by fixing $\lambda_1 = 1$, $\lambda_2 = 1 - \Delta$, and $\lambda_n = 0$, with the interior eigenvalues drawn uniformly at random in descending order. Table \ref{tab:tabsota} summarizes the experimental results, evaluated in terms of both computation time and the number of matrix-vector products.
	
	\textbf{Accelerated Efficiency over Power Iterations.} 
	The proposed Split-Merge algorithm exhibits a consistent advantage over power iteration and its variants. Relative to the basic power method, Split-Merge achieves an overall speed-up factor ranging from $4.59\times$ to $9.16\times$, with the relative efficiency increasing as the problem difficulty scales (i.e., at $n=4096, \Delta=10^{-3}$). Furthermore, Split-Merge circumvents the hyperparameter sensitivity inherent to momentum-based methods. While the Power+M algorithm necessitates the exact optimal parameter $\beta^*$, a slight perturbation to $\beta = 0.9\beta^*$ degrades its performance. For example, at $n=4096$ and $\Delta=10^{-3}$, the runtime of Power+M (near-ideal) is $7.92 \times 10^{-1}$s, whereas Split-Merge requires only $2.58 \times 10^{-1}$s. This execution time yields a roughly $3\times$ speed-up over the near-ideal case and surpasses even the idealized Power+M baseline ($2.67 \times 10^{-1}$s), which excludes the computational overhead of estimating the second-largest eigenvalue.
	
	\textbf{Iteration Complexity vs. Computational Cost.} 
	The inclusion of MV counts provides critical insight into the per-iteration complexity of subspace methods. As expected, subspace methods such as Lanczos and LOBPCG require the fewest MVs to converge (e.g., $12$--$37$ MVs and $14$--$45$ MVs, respectively). However, a lower MV count does not strictly translate to faster computation times. At $n=1024$ and $\Delta=10^{-1}$, Split-Merge requires slightly more MVs (17) than LOBPCG (14), yet it achieves a computation time of $8.28 \times 10^{-4}$s compared to LOBPCG's $6.85 \times 10^{-3}$s, yielding a speed-up of more than $8.2\times$. This discrepancy reflects a structural distinction of Split-Merge: by avoiding the Rayleigh-Ritz projections and explicit orthogonalization steps required by LOBPCG and JD, it maintains a lower per-iteration overhead. This theoretical advantage translates directly into substantial empirical gains; for instance, at $n=2048$ and $\Delta=10^{-2}$, Split-Merge operates approximately $5.2\times$ faster than the JD method. Although the classical Lanczos method requires the least computation time in most test scenarios, Split-Merge remains highly competitive, often securing the first or second-best time metrics across various setups.
	
	\textbf{Comparisons with Difference-Based Approaches.} 
	Split-Merge significantly outperforms difference-based approaches. At $n=2048$ and $\Delta=10^{-2}$, for instance, Split-Merge solves the problem in $7.68 \times 10^{-3}$s, realizing speed-ups of approximately $11.7\times$ and $3.9\times$ over the Newton and L-BFGS methods, respectively.
	
	\textbf{Comparisons with Standard Software Baselines.}
	To ensure a rigorous evaluation against industry-standard routines, we benchmarked our approach against MATLAB's native \texttt{eigs}. Despite the inherent execution advantages of its highly optimized, compiled architecture, the Split-Merge algorithm demonstrates competitive performance. For instance, in moderately sized problems (e.g., $n=1024, \Delta=10^{-1}$), Split-Merge ($8.28 \times 10^{-4}$s) computes the solution nearly $2.6\times$ faster than the \texttt{eigs} routine ($2.19 \times 10^{-3}$s).
	Collectively, these comprehensive metrics firmly establish Split-Merge as a highly efficient and hyperparameter-free alternative for practical eigenvalue computations.
	
	\textbf{Impact of Diagonal Shifts.}
	For indefinite matrices, a diagonal shift $\boldsymbol{A} + \eta \boldsymbol{I}$ enforces positive semi-definiteness. Using an indefinite configuration ($n = 2048$, $\Delta = 10^{-1}$, $\lambda_n = -10$), we vary $\eta$ from $-\lambda_n$ to the Gershgorin bound \cite{golub2013matrix}. As shown in \cref{fig:shift_eta_impact}, the iteration complexity of the power method grows linearly with $\eta$ due to the degradation of its static spectral ratio $(\lambda_i + \eta)/(\lambda_1 + \eta)$ toward $1$. Conversely, Split-Merge maintains a highly consistent iteration complexity, matching the Lanczos baseline. Corroborating \cref{rem:shift_insensitivity}, the dynamic root $\tau_k$ adapts to the shifted spectrum, enabling $\eta$ to effectively cancel within the attenuation factor. This ensures a convergence rate largely insensitive to the shift magnitude.
	}
	
	\begin{figure}[t]
		\centering
		\begin{subfigure}[b]{0.32\linewidth}
			\centering
			\includegraphics[width=\textwidth]{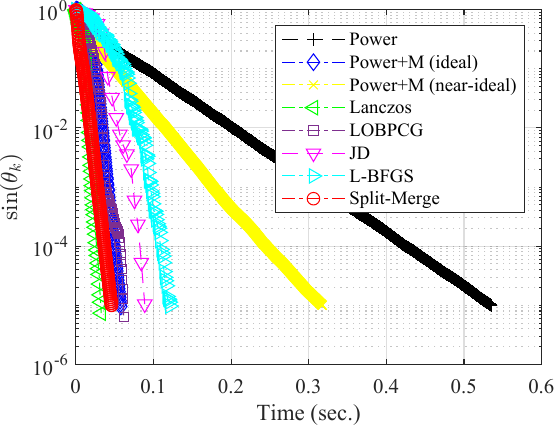}
			\caption{Kuu}
			\label{fig:sin_iter_kuu}
		\end{subfigure}
		\hfill
		\begin{subfigure}[b]{0.32\linewidth}
			\centering
			\includegraphics[width=\textwidth]{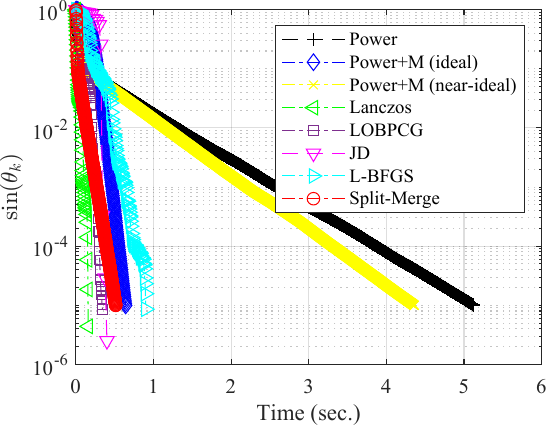}
			\caption{Andrews}
			\label{fig:sin_time_Andrews}
		\end{subfigure}
		\hfill
		\begin{subfigure}[b]{0.32\linewidth}
			\centering
			\includegraphics[width=\textwidth]{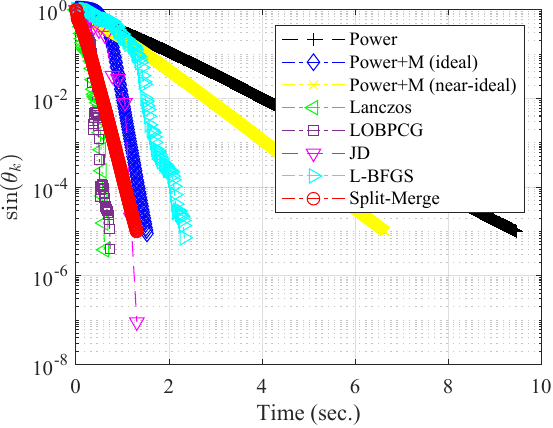}
			\caption{boneS01}
			\label{fig:sin_time_boneS01}
		\end{subfigure}
		\hfill
		\begin{subfigure}[b]{0.32\linewidth}
			\centering
			\includegraphics[width=\textwidth]{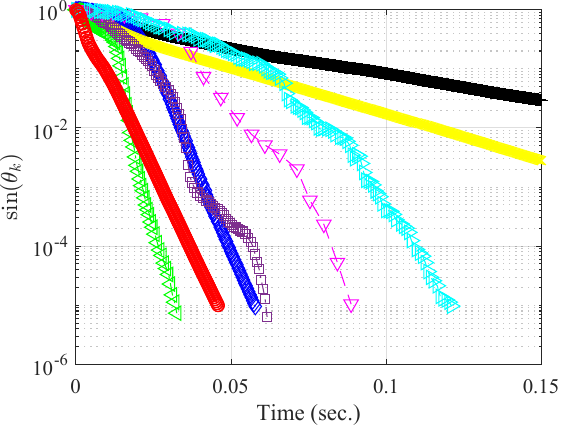}
			\caption{Kuu (Zoomed in)}
			\label{fig:sin_time_kuu_zoom}
		\end{subfigure}
		\hfill
		\begin{subfigure}[b]{0.32\linewidth}
			\centering
			\includegraphics[width=\textwidth]{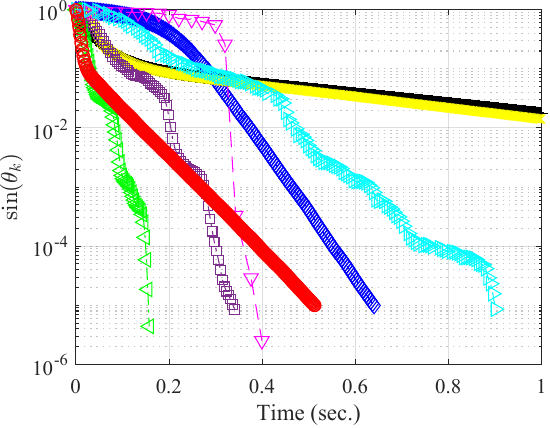}
			\caption{Andrews (Zoomed in)}
			\label{fig:sin_time_Andrews_zoom}
		\end{subfigure}
		\hfill
		\begin{subfigure}[b]{0.32\linewidth}
			\centering
			\includegraphics[width=\textwidth]{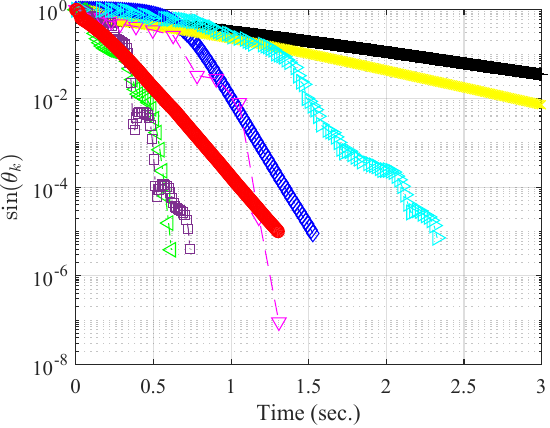}
			\caption{boneS01 (Zoomed in)}
			\label{fig:sin_time_boneS01_zoom}
		\end{subfigure}
		\caption{\tp \textbf{Comparison of different methods on three SuiteSparse benchmark matrix datasets}. The speed-up achieved by Split-Merge method is highlighted for each dataset: (a) Kuu: \textbf{11.67}, (b) Andrews: \textbf{9.98}, and (c) boneS01: \textbf{7.26}.}
		\label{fig:sin_time_suitesparse}
	\end{figure}
	
	\subsection{SuiteSparse Matrix Datasets}

	{
	\tp
	We further evaluate the empirical performance of our algorithm using three large-scale PSD benchmark matrices from the SuiteSparse Matrix Collection\footnote{\tp Available at \url{https://sparse.tamu.edu/}} \cite{suitesparse}. The selected matrices vary significantly in dimension: \texttt{Kuu} ($n=7,102$), \texttt{Andrews} ($n=60,000$), and \texttt{boneS01} ($n=127,224$). To establish a rigorous ground truth, we utilize MATLAB's built-in \texttt{eigs} function to compute the dominant eigenpair for each matrix.
	
	\cref{fig:sin_time_suitesparse} reports the approximation error $\sin(\theta_k)$ over computation time. The trajectories show that the time-efficiency of Split-Merge is comparable to established subspace methods (Lanczos, LOBPCG, and JD), while outperforming both Power+M and L-BFGS. Additionally, the early-stage plots (\cref{fig:sin_time_kuu_zoom}) indicate that Split-Merge attains low-to-medium precision regimes more rapidly.
	
	Finally, when compared to the standard power method, Split-Merge yields substantial computational savings across all test cases. Notably, it achieves an impressive speed-up of over $10\times$ on the \texttt{Kuu} matrix. This robust acceleration scales seamlessly to the larger datasets, consistently delivering speed-ups exceeding $7\times$ on both the \texttt{Andrews} and \texttt{boneS01} matrices.
	}
	
	\begin{table}[tbp]
	\centering
	\caption{\tp \textbf{Comparison of average computation time (sec.) and the number of matrix-vector products across three real-world datasets}. Bold and underlined values indicate the \textbf{best} and \underline{second-best} results, respectively. The proposed method is highlighted as \colorbox{smblue}{Split-Merge}. Additionally, the \textit{speed-up} factor of the Split-Merge algorithm relative to the power method is reported.}
	\label{tab:ucl_benchmark}
	\renewcommand{\arraystretch}{1.12}
	\setlength{\tabcolsep}{4pt}
	\footnotesize
	\begin{tabular}{l|cc|cc|cc}
		\toprule
		\multirow{2}{*}{\textbf{Method}} & \multicolumn{2}{c}{{\texttt{Gisette}}} & \multicolumn{2}{c}{{\texttt{Arcene}}} & \multicolumn{2}{c}{{\texttt{GeneExp}}} \\
		\cmidrule(lr){2-3}\cmidrule(lr){4-5}\cmidrule(lr){6-7}
		& \textbf{Time} & \textbf{MVs} & \textbf{Time} & \textbf{MVs} & \textbf{Time} & \textbf{MVs} \\
		\midrule
		\textbf{Power} & 0.36 & 86 & 0.52 & 31 & 6.98 & 98 \\
		\textbf{Power$+$M (ideal)} & 0.20 & 47 & 0.44 & 25 & 3.11 & 43 \\
		\textbf{Power$+$M (near-ideal)} & 0.20 & 46 & 0.48 & 27 & 3.77 & 52 \\
		\textbf{Lanczos} & \textbf{0.04} & \textbf{9} & \textbf{0.11} & \textbf{6} & \textbf{0.73} & \textbf{10} \\
		\textbf{LOBPCG} & 0.09 & \underline{13} & 0.22 & \underline{8} & 1.45 & \underline{13} \\
		\textbf{JD} & 0.13 & 20 & 0.21 & 10 & 1.49 & 18 \\
		\textbf{L-BFGS} & 0.34 & 107 & 0.92 & 83 & 5.41 & 106 \\
		\rowcolor{smblue}\textbf{Split-Merge} & \underline{0.08} & 18 & \underline{0.17} & 10 & \underline{1.36} & 19 \\
		\midrule
		\texttt{eigs} & 0.12 & -- & 0.45 & -- & 1.77 & -- \\
		\midrule
		\textit{Speed-up} & \multicolumn{2}{c}{\textit{4.79$\times$}} & \multicolumn{2}{c}{\textit{3.10$\times$}} & \multicolumn{2}{c}{\textit{5.13$\times$}} \\
		\bottomrule
	\end{tabular}
\end{table}
	
	\subsection{Real-World Datasets}
	
	To evaluate the practical efficacy of the Split-Merge algorithm in real-world PCA tasks, we apply it to three high-dimensional datasets from the UCI Machine Learning Repository\footnote{Available at \url{https://archive.ics.uci.edu/}}: \texttt{Gisette} \cite{gisette_170}, \texttt{Arcene} \cite{arcene_167}, and \texttt{GeneExp} \cite{gene_expression_cancer_rna-seq_401}. Prior to evaluation, the features of each dataset are standardized to zero mean and unit variance.
	
	\textbf{Datasets Overview.}
	\begin{itemize}
		\item \texttt{Gisette} \cite{gisette_170}: Consists of 13,500 samples with 5,000 features, originally constructed for handwritten digit recognition.
		\item \texttt{Arcene} \cite{arcene_167}: Consists of 900 samples with 10,000 features, designed to classify cancerous versus normal patterns in mass-spectrometric data.
		\item \texttt{GeneExp} \cite{gene_expression_cancer_rna-seq_401}: Contains 801 samples with 20,531 features, focusing on cancer classification through gene expression analysis.
	\end{itemize}
	
	{
	\tp
	\cref{tab:ucl_benchmark} summarizes the average computation time and matrix-vector products evaluated over 500 independent trials. The Split-Merge algorithm consistently secures the second-best computational time across all datasets, surpassed only by the Lanczos baseline, while maintaining a performance advantage over both momentum-accelerated power variants and L-BFGS. Crucially, these metrics empirically validate the core theoretical mechanism: replacing the trivial assignment $\boldsymbol{v}_k = \boldsymbol{0}$ (which recovers the standard power method) with the exact orthogonal projection $\boldsymbol{v}_k$ defined in \cref{eq:v_sm} overcomes the convergence bottleneck. This principled modification directly translates into massive computational savings, delivering up to a $5.13\times$ speed-up over the classical scheme in high-dimensional machine learning scenarios.
	}
	
	\section{Conclusion}\label{sec:conclusion}
	
	{
	\tp
	This paper recasts the classical power method within a first-order optimization framework. By extending the analysis from the equivalence at the fixed stepsize $\alpha = 1/2$ to the entire interval $\alpha \in (0, 1)$, we establish almost sure convergence to the global minimizer and characterize the local linear rate. This perspective not only exposes the power method as a conservative instance but also quantifies its asymptotic sub-optimality.
	To overcome the isotropic limitations of classical first-order updates, we develop the Split-Merge algorithm under the majorization-minimization principle.
	By implicitly splitting the matrix to construct a curvature-aware surrogate and maximizing the descent via a tractable relaxation, the algorithm extracts auxiliary vectors that merge the factors, returning a matrix-free iteration.
	Characterizing this process as an adaptive spectral filter reveals a spectral peeling mechanism that suppresses localized eigenspaces, breaking the linear convergence bottleneck.
	Supported by evaluations across different datasets, Split-Merge delivers computational efficiency comparable to subspace methods.
	Subsequent research will extend this approach to generalized eigenvalue problems and the simultaneous extraction of multiple leading eigenvectors on matrix manifolds.
	}
	
	\section*{Acknowledgments}
	{\tp The authors are grateful to the associate editor and the two anonymous referees for their very helpful comments and suggestions.}
	
	\bibliographystyle{siamplain}
	\bibliography{references}
\end{document}